\documentclass[lineno]{biometrika}

\usepackage{times}
\usepackage{natbib}
\usepackage{amssymb}
\usepackage{amsmath}
\usepackage{caption}

\newcommand{\biblist}{
\bibliographystyle{biometrika}
\setlength{\bibsep}{0cm}
}

\usepackage[plain,noend]{algorithm2e}
\usepackage{xcolor}

\makeatletter
\renewcommand{\algocf@captiontext}[2]{#1\algocf@typo. \AlCapFnt{}#2} 
\renewcommand{\AlTitleFnt}[1]{#1\unskip}
\def\@algocf@capt@plain{top}
\renewcommand{\algocf@makecaption}[2]{%
  \addtolength{\hsize}{\algomargin}%
  \sbox\@tempboxa{\algocf@captiontext{#1}{#2}}%
  \ifdim\wd\@tempboxa >\hsize
    \hskip .5\algomargin%
    \parbox[t]{\hsize}{\algocf@captiontext{#1}{#2}}
  \else%
    \global\@minipagefalse%
    \hbox to\hsize{\box\@tempboxa}
  \fi%
  \addtolength{\hsize}{-\algomargin}%
}
\makeatother

\newcommand{\Cov}{\mathrm{cov}}
\newcommand{\E}{E}

\newcommand{\N}{{\mathbb N}}
\newcommand{\R}{{\mathbb{R}}}
\newcommand{\Z}{{\mathbb Z}}
\newcommand{\spCov}{PP^\T + \eta_1^2 I_k}
\newcommand{\tCov}{V}
\newcommand{\spCC}{P q}

\newcommand{\PLS}{partial least squares}
\newcommand{\txtCG}{conjugate gradient }

\newcommand{\ARMA}{autoregressive moving average process}

\def\T{{ \mathrm{\scriptscriptstyle T} }}
\def \calL {\mathcal L}

\begin{document}
\nolinenumbers


\markboth{M. Singer, T. Krivobokova, A. Munk \and B. de Groot}{Partial least squares for dependent data}

\title{Partial least squares for dependent data}

\author{MARCO SINGER, TATYANA KRIVOBOKOVA, AXEL MUNK}
\affil{Institute for Mathematical Stochastics, Georg-August-Universit\"at G\"ottingen, Goldschmidtstr. 7, 37077 G\"ottingen, Germany
\email{msinger@gwdg.de} 
\email{tkrivob@uni-goettingen.de}
\email{munk@math.uni-goettingen.de}}
\author{\and BERT DE GROOT}
\affil{Max Planck Institute for Biophysical Chemistry, Am Fassberg 11, 37077 G\"ottingen, Germany 
\email{bgroot@gwdg.de}}

\maketitle

\begin{abstract}
We consider the partial least squares algorithm for dependent data and study the consequences of ignoring the dependence both theoretically and numerically. Ignoring non-stationary dependence structures can lead to inconsistent
estimation, but a simple modification leads to consistent estimation. A protein dynamics example illustrates the superior predictive
power of the method.
\end{abstract}


\pagenumbering{arabic}

\section{Introduction}
\label{sec:Introduction}
The partial least squares algorithm introduced by 
\citet{Wold66} is a powerful regularized regression tool. It is an iterative technique, which is,
unlike most similar methods, non-linear in the response variable. Consider a linear regression model
\begin{equation}
	\label{eq:linMod}
y = X \beta + \varepsilon,
\end{equation}
where $y\in\mathbb{R}^n$ , $X\in\mathbb{R}^{n\times k}$, $\beta \in \R^k$ and the error term $\varepsilon\in \mathbb{R}^n$ is a vector of $n$ independent and identically distributed random variables. To estimate the unknown coefficients $\beta$ with partial least squares, a base of $i\leq k$ weight vectors
$\widehat{w}_1,\ldots,\widehat{w}_i$ is iteratively
constructed. First, the data are centered, i.e., $y$ and the columns of $X$ are transformed to have mean zero. Then the  first vector $\widehat{w}_1$ is obtained by maximizing the empirical
covariance between $X w$ and $y$ in $w \in \R^k$, subject to $\|w\|=1$. Afterwards,
the data are projected into the space orthogonal to $X
\widehat{w}_1$ and the procedure is iterated. The $i$th \PLS{} estimator
$\widehat{\beta}_i$ for $\beta$ is obtained by performing a least
squares regression of $y$ on $X$, constrained to the subspace spanned by the columns of
$\widehat{W}_i=(\widehat{w}_1,\ldots,\widehat{w}_i)$. \citet{Hel88} summarizes the \PLS{} iterations in
two steps via
\begin{eqnarray}
	\label{eq:PLSReg}
        \widehat{w}_{i+1} & = &b - A \widehat{\beta}_i,\;\;\widehat{\beta}_0=0,
	 \\
	\notag
	\widehat{\beta}_i & =& \widehat{W}_i(\widehat{W}_i^\T A \widehat{W}_i)^{-1} \widehat{W}_i^\T b,
\end{eqnarray}
with $b=n^{-1} X^\T y$ and $A=n^{-1}X^\T X$, under the assumption that $(\widehat{W}_i^\T A \widehat{W}_i)^{-1}$ exists. The regularisation is achieved by early stopping, that is, by taking
$i\leq k$.

Alternatively, $\widehat{\beta}_i$ can be defined using the
fact that $\widehat{w}_i\in{\mathcal{K}}_i(A,b)$, where
${\mathcal{K}}_i(A,b)$ is a Krylov space, that is, a space
spanned by $\left\{A^{j-1}b\right\}_{j=1}^i$ \citep[][]{Hel88}. Then, one can define \PLS{}
estimators as
$
\widehat{\beta}_i =\arg\min_{\beta\in {\mathcal{K}}_i(A,b)}(y-X\beta)^\T(y-X\beta).
$
There is also a direct correspondence between partial least squares and the 
conjugate gradient method with early stopping for the solution of
$A\beta=b$. 

\citet{FrankFriedman} and \citet{Far05} find the \PLS{} algorithm to be competitive with regularized regression
techniques, such as principal component regression, lasso or ridge
regression, in terms of the mean squared prediction
error. Also, the optimal number of partial least squares base components is often much lower than that of principal components regression, as found in \citet{Almoy96}.

Partial least squares regression has a long and successful history in various application areas, see e.g.,  
\citet{Hul99}, \citet{Lob01}, \citet{Ngu02}. However, the statistical
properties of this algorithm have been little studied, perhaps because of the non-linearity of \PLS{} estimators in the response
variable.
Some attempts to understand properties of \PLS{} can be found
 in \citet{Hoes88}, \citet{Pha02} and \citet{Krae07}. Their almost sure convergence was established by \citet{Naik00}.
 For kernel partial least squares, \citet{Blan10a} obtained convergence in probability results by early stopping. For the closely linked kernel conjugate gradient algorithm, \citet{Blan10b} established order-optimal convergence rates dependent on the regularity of the target function.
\citet{Del12} compared theoretically the population and sample properties of the \PLS{}
algorithm for functional data. 

Regression techniques typically assume independence of responses, but this is often violated, for example, if the data are observed over
time or at dependent spatial locations. 
We are not aware of any treatment of the \PLS{} algorithm for dependent observations. 
In this work we propose a modification of \PLS{} to
deal with dependent observations and study the theoretical
properties of \PLS{} estimators under general dependence in the data. In particular, we 
quantify the influence of ignored dependence. 

Throughout the paper we denote by $\|\cdot\|_\calL$ the spectral and by $\|\cdot\|$ the Frobenius norm for matrices, $\|\cdot\|$ also denotes the Euclidean norm for vectors.

\section{\PLS{} under dependence}
\label{sec:PLS}
\subsection{Latent variable model}
In many applications the standard linear model (\ref{eq:linMod}) is too
restrictive. For example, if a covariate that is relevant for the response cannot be observed or measured directly, so-called  latent variable or structural equation models are considered \citep[][]{Skro06}: it is assumed that
$X$ and $y$ are linked by $l\leq k$ latent vectors and the
remaining vectors in the $k$-dimensional column space of $X$ do not contribute to
$y$. This can be interpreted as if the latent
components are of interest, but only $X$, which contains some unknown nuisance information, can be
measured. Such models are relevant in modelling of chemical
  \citep{Wol01}, economic \citep{Hahn02} and social data \citep{Gold72}.

We consider a latent variable model with the covariates $X$ and response
$y$ connected via a matrix of latent variables $N$,
\begin{equation}
\begin{aligned}
	X &= \tCov (N P^\T + \eta_1F), \\ 
	y &=\tCov (N q + \eta_2f),
\end{aligned}
\label{eq:genModel}
\end{equation}

where $N$
and
$F $
are an $n \times l$-dimensional and an $n \times k$-dimensional random matrix, respectively, and $f$ is an $n$-dimensional random vector. The random elements $N$,
$F$, $f$ can have different
distributions, but are independent of each other, with all entries being
independent and identically distributed with expectation zero and unit variance. The matrix
 $P \in \R^{k \times l}$ and vector $q\in \R^l$ are
deterministic and unknown, along with the real-valued parameters $\eta_1,\eta_2\geq 0$. We assume that $n\geq k \geq l$ and that $\mathrm{rank}(N) = \mathrm{rank}(P) = l$, $\mathrm{rank}(F) = k$ almost surely. 

The matrix $\tCov \in\R^{n
  \times n}$ is a deterministic symmetric matrix, such that
$\tCov^2$ is a positive definite covariance matrix. 
If $\tCov\neq I_n$, then $X$ in model (\ref{eq:genModel}) can be seen as the matrix form of a $k$-dimensional time series
$\{X_t\}_{t=1}^n$ ($X_t\in\mathbb{R}^k$) and
$y$ can be seen as a real-valued time series
$\{y_t\}_{t=1}^n$. The covariance matrix $\tCov^2$ determines the
dependence between
observations, which might be non-stationary. {We will call $\tCov^2$ the temporal covariance matrix of $X$ and define $\Sigma^2 = \spCov$.} Setting $l=k$, $\eta_1=0$ reduces model (\ref{eq:genModel}) to standard linear regression with dependent observations. 

The latent variables $N$ connect $X$ to $y$, whereas $F$ can
be considered as noise, thus giving a model where not all directions
in the column space of $X$ are important for the prediction of
$y$. The representation (\ref{eq:genModel}) highlights practical
settings where the \PLS{} algorithm is expected to outperform principal
component regression and
similar techniques. In particular, if the covariance of $\eta_1 F$
dominates that of $NP^\T$, then the first principal components will be largely uncorrelated to $y$. In contrast,
the first \PLS{} basis components should by definition be able to recover relevant latent components.

The \PLS{} algorithm is run as described in Section \ref{sec:Introduction} with matrix
$X$ and vector $y$ defined in (\ref{eq:genModel}). 
If $\eta_1 = 0$, then model (\ref{eq:linMod}) is correctly specified
with {$q=P^\T \beta$} and the \PLS{} estimator (\ref{eq:PLSReg})
estimates $\beta$. If $\eta_1>0$, then model (\ref{eq:linMod}) is
misspecified and  
$\beta(\eta_1)=\Sigma^{-2} \spCC$ is rather estimated. Note that $\beta(0)=\beta$.

In the standard \PLS{} algorithm it is assumed that $\tCov=I_n$. In the 
subsequent sections we aim to quantify the influence of $\tCov\neq I_n$, which
is ignored in the algorithm.


\subsection{Population and sample partial least squares}
\label{subsec:samplePLS}
The population \PLS{} algorithm for independent observations was first introduced by
\citet{Hel90}. Under model (\ref{eq:genModel}) we modify the definition of the population \PLS{} basis vectors as
$$
	w_i = \arg \max\limits_{\substack {w \in \R^k\\
            \|w\|=1}} \frac{1}{n^2} \sum\limits_{t,s=1}^n
        \Cov(y_t-X_t^\T\beta_{i-1} ,  X_s^\T
        w),\;\;\beta_0=0,
$$
where $\beta_i\in \R^k$ are the population \PLS{} regression
coefficients. The average covariances over observations are taken,
since the data are neither
independent nor identically distributed if $\tCov^2\neq I_n$. 
Solving this optimization problem implies that the basis
vectors $w_1,\ldots,w_i$ span the Krylov space $\mathcal{K}_i(\Sigma^2,\spCC)$:
see the Supplementary Material. In particular, under model (\ref{eq:genModel}), the Krylov
space in the population turns out to
be independent of the temporal covariance $\tCov^2$ for all $n \in \N$.

For a given Krylov space, the population \PLS{} coefficients are obtained as
\begin{equation*}
	\beta_i=\arg\min\limits_{\beta \in \mathcal{K}_i(\Sigma^2 , \spCC)} \E \left\{\frac{1}{n}\sum\limits_{t=1}^n   \left(
		y_t - X_t^\T \beta
	\right)^2\right\}. 
\end{equation*}%
It is easy to see that the solution to this problem is 
$$
	\beta_i = K_i \left(K_i^\T \Sigma^2 K_i \right)^{-1}
        K_i^\T \spCC,\;\;K_i=(\spCC,\Sigma^2 \spCC,\ldots,\Sigma^{2(i-1)}\spCC),
$$
which is independent of $\tCov^2$ for all $n \in \N$.

To obtain the sample \PLS{} estimators $\widehat\beta_i$,  $\Sigma^2$ and $Pq$ are replaced by
estimators. In the standard \PLS{} algorithm, under independence of
observations, $\Sigma^2$ and $Pq$ are estimated by unbiased
estimators $n^{-1}X^\T X$
and $n^{-1}X^\T y$, respectively. However, if the observations are
dependent, such naive estimators can lead to $L_2$-inconsistent estimation,
as the following theorem shows.
\begin{theorem}
\label{th:FirstComponentError}
	Let the model (\ref{eq:genModel}) hold and the fourth moments of $N_{1,1}$, $F_{1,1}$ exist. Define
$A=\|\tCov\|^{-2}X^\T X,\;\;\;b=\|\tCov\|^{-2}X^\T y$. Then
\begin{eqnarray*}
\E\left(\left\|\Sigma^2 - A \right\|^2\right)& =&\frac{\|\tCov^2\|^2}{\|\tCov\|^4}\left(C_{A}+\sum\limits_{t=1}^n\frac{\|\tCov_t\|^4}{\|\tCov^2\|^2} c_{A}\right)\\
	\E\left(\|\spCC-b\|^2\right) &=&
	\frac{\|\tCov^2\|^2}{\|\tCov\|^4}\left(C_{b}+\sum\limits_{t=1}^n\frac{\|\tCov_t\|^4}{\|\tCov^2\|^2} c_{b}\right),
\end{eqnarray*}
where 
\begin{eqnarray*}
C_{A}&=&
			\|P\|^4 + \|P^\T P\|^2 + 4 \eta_1^2 \|P\|^2 + 
			\eta_1^4 k (1+k)\\
c_{A}&=&
		 	\left\{
		 		\E\left(N_{1,1}^4\right)-3
		 	\right\}\sum\limits_{i=1}^l \|P_i\|^4  + \left\{
		 		\E\left(F_{1,1}^4\right) -3
		 	\right\}\eta_1^4 k \\
C_{b}&=&\|\spCC\|^2 + \|P\|^2\|q\|^2 + \eta_1^2 k \|q\|^2 + \eta_1^2 \eta_2^2 k 
+ \eta_2^2 \|P\|^2\\ 
c_{b}&=&	\left\{
				\E\left( N_{1,1}^4
				\right) -3
		\right\}	\sum\limits_{i=1}^l \|P_i\|^2 q_i^2 
			\end{eqnarray*}
and $V_t$ denotes the $t$-th column of matrix $V$.
\end{theorem}

The scaling factors in $A$ and $b$ have no influence on the
sample \PLS{} estimators in (\ref{eq:PLSReg}), so that replacing $n^{-1}$
with $\|\tCov\|^{-2}$ does not affect the algorithm and both $A$ and $b$ are unbiased estimators for $\Sigma^2$ and $P q$, respectively.

If $\E(N_{1,1}^4) =  \E(F_{1,1}^4) = 3$, then constants $c_{A}$ and $c_{b}$ vanish, simplifying expressions for the mean squared error of $A$ and $b$. This is satisfied, for example, for the standard
  normal distribution. Thus, these terms can be interpreted as a penalization for non-normality. 

Finally, $
	\sum_{t=1}^n \|\tCov_t\|^4 \leq 
	\sum_{t,s=1}^n \left(\tCov_t^\T \tCov_s\right)^2 = \left\|\tCov^2\right\|^2
$ implies that the convergence rate of both estimators is driven by the ratio of Frobenius norms $\|\tCov\|^{-2} \|\tCov^2\|$. In particular, if $\|\tCov\|^{-2} \|\tCov^2\|$ converges to
zero, then the elements of the population Krylov space $\Sigma^2$ and $Pq$ can be estimated consistently. This
is the case, for example, for independent observations with $\tCov =
I_n$, since $\|I_n^2\| = \|I_n\| =
n^{1/2}$. However, if $\|\tCov\|^{-2} \|\tCov^2\|$ fails to converge to zero, ignoring the
temporal dependence $\tCov^2$  may lead to inconsistent
estimation. 

\section{Properties of \PLS{} estimators under
  dependence}
\label{sec:theory}
\subsection{Concentration inequality for \PLS{} estimators}
\label{sec:Concentration}
In this section we apply techniques of \citet{Blan10b}, who derived
convergence rates of the kernel \txtCG algorithm, which
is closely related to kernel \PLS{}. Both algorithms approximate the
solution on Krylov subspaces, but employ different norms. In particular, \citet{Blan10b} have shown that 
if the conjugate gradient algorithm is stopped early, the convergence in
probability of the kernel \txtCG estimator to the true regression
function can be obtained for bounded kernels. Moreover, the 
convergence is order-optimal, depending on the regularity of the
target function. These results hold for independent identically
distributed observations.

We avoid the nonparametric setting of \citet{Blan10b} and
study a standard linear \PLS{} algorithm with a fixed dimension $k$ of the
regression space. We allow the observations to be dependent,
and, instead of a bounded kernel, consider unbounded random variables
with moment conditions. In this setting we derive concentration
inequalities for \PLS{} estimators that allow us to quantify the influence
of the temporal covariance.

Regularization of the \PLS{} solution is achieved by early stopping,
which is characterized by the discrepancy principle, i.e., we stop at the first index
$0<a_0\leq a$ such that 
\begin{equation}
	\label{eq:discrepancy0}
\left\|A^{1/2} \widehat{\beta}_{a_0} -A^{-1/2}b\right\| \leq \tau(\delta \|\widehat{\beta}_{a_0}\| + \epsilon),
\end{equation}
for $\delta,\epsilon >0$ defined in Theorem
\ref{th:ExpBlanchard}, and some $\tau \geq 1$. Here $a$ denotes the maximal dimension of the sample Krylov space
$\mathcal{K}_i(A,b)$ and almost surely equals $a= l+(k-l)\mathbb{I}(\eta_1>0)$, where $\mathbb{I}(\cdot)$ denotes an indicator function. For technical reasons we stop at $a^*=a_0-1$ if $p_{a_0}(0)
\geq \zeta \delta^{-1}$, where $p_i$ is a polynomial of degree $i-1$
with $p_i(A)b = \widehat{\beta}_i$ and $\zeta < \tau^{-1}$.
The existence of such polynomials was proved by \citet{Pha02}. If (\ref{eq:discrepancy0}) never holds, $a^*= a$ is
taken. With this stopping index we get the following concentration inequality.
\begin{theorem}
\label{th:ExpBlanchard}
Assume that model (\ref{eq:genModel}) with $\eta_1>0$ holds and that the fourth moments of $N_{1,1}$, $F_{1,1}$ exist. Furthermore, $a^*$ satisfies (\ref{eq:discrepancy0}) with $\tau\geq 1$, $\zeta<\tau^{-1}$. For $\nu \in (0,1]$ let $\delta=\nu^{-1/2}\|\tCov\|^{-2} \|\tCov^2\| C_\delta$ and $\epsilon=\nu^{-1/2}\|\tCov\|^{-2} \|\tCov^2\| C_\epsilon$, such that $\delta,\epsilon\rightarrow 0$, where
\begin{equation*}
	C_\delta= \left(
		2C_A+2c_A\right)^{1/2},	\;\;
	C_\epsilon= \left(
		2C_b+2c_b\right)^{1/2},
\end{equation*}
with $C_A$, $c_A$, $C_b$ and $c_b$ given in Theorem \ref{th:FirstComponentError}. Then with a probability at least $1- \nu$,
\begin{equation}
	\label{eq:BlanchardInequality}
		\left\|
			\widehat{\beta}_{a^*} - \beta(\eta_1)
		\right\| \leq \frac{\|\tCov^2\|}{\|\tCov\|^2}\left\{ c_1(\nu) + \frac{\|\tCov^2\|}{\|\tCov\|^2} c_2(\nu) \right\},
\end{equation}
where 
\begin{align*}
	c_1(\nu) &= \nu^{-1/2} c(\tau,\zeta)\|\Sigma^{-1}\|_\calL \left(
 		C_\epsilon + \|\Sigma\|_\calL\|\Sigma^{-3} P q\| C_\delta 
 	\right)\\
 	c_2(\nu) & = \nu^{-1}c(\tau,\zeta)\|\Sigma^{-1}\|_\calL \left(
 	C_\epsilon C_\delta + \|\Sigma^{-3} P q \| C_\delta^2
 	\right),
\end{align*}
for some constant $c(\tau,\zeta)$ that asymptotically depends only on $\tau$ and $\zeta$.
\end{theorem}

If $N_{1,1}$,$F_{1,1}$, $f_1\sim{\mathcal{N}}(0,1)$, then the expressions for $C_\delta$ and $C_\epsilon$ are simplified and  
  the scaling factor of $c_1(\nu)$ and $c_2(\nu)$ can be improved from
  $\nu^{-1/2}$ to $\log(2/\nu)$, which is achieved by
  using an exponential inequality proved in Theorem 3.3.4 of \citet{Yurinski95}. 

Theorem \ref{th:ExpBlanchard} states that the convergence rate of the optimally stopped \PLS{} estimator $\widehat\beta_{a^*}$
to the true parameter
$\beta(\eta_1)$ is driven by the ratio of
the Frobenius norms of $\tCov^2$ and $\tCov$, similar to the results of
Theorem \ref{th:FirstComponentError}.
In particular, if the data are independent with $\tCov=I_n$ then $\widehat{\beta}_{a^*}$
is square-root consistent. In this case $c_2(\nu)$ is asymptotically negligible. Note that the theorem excludes the case that $\|V\|^{-2}\|V^2\|$ does not converge to $0$.

\subsection{Properties of $\widehat\beta_1$ under dependence}
Non-linearity in the response variable of $\widehat\beta_i$ hinders its
standard statistical analysis, as no closed-form expression for the mean square error of $\widehat\beta_i$ is available and concentration
inequalities similar to (\ref{eq:BlanchardInequality}) are the only
results on the convergence rates of \PLS{} estimators, to the best of
our knowledge. However, if the ratio of $\|\tCov^2\|$ and $\|\tCov\|^2$ does not
converge to zero, Theorem \ref{th:ExpBlanchard} does not hold.

In this section we study the first \PLS{} estimator $\widehat\beta_1$, for
several reasons. First, the explicit expression for its mean square error can be derived. Second, if there is
only one latent component that links $X$ and $y$, i.e., $l=1$ in (\ref{eq:genModel}), then consistent estimation of $\beta_1$ is
crucial. Finally, $\widehat\beta_1$ is collinear to the direction of the
maximal covariance between $X$ and $y$ given by $\widehat{w}_1$,
which is important for the interpretation of the \PLS{} model in
applications, see \citet{KrivGroot12}.
The next theorem gives conditions under which $\widehat\beta_1$
is an inconsistent estimator of $\beta_1$.
\begin{theorem}
\label{th:Inconsistency}
Assume that model (\ref{eq:genModel}) holds, $k>1$ and eighth moments of $N_{1,1}$, $F_{1,1}$, $f_1$ exist.
Furthermore, suppose that the ratio $\|\tCov\|^{-2}\|\tCov^2\| $ does not converge to zero as $n\rightarrow\infty$. Then, for either $l>1$, $\eta_1 \geq 0$ or $l=1$, $\eta_1>0$, $\widehat{\beta}_1$ is an inconsistent estimator for $\beta_1$.
\end{theorem}

The case $l=1$, $\eta_1=0$ not treated in Theorem
  \ref{th:Inconsistency} corresponds to the standard
  linear regression model with a single covariate, so the \PLS{}
  estimator coincides with the ordinary least squares estimator, see \citet{Hel88}.
  
Hence, if there is only one latent component in the model, i.e., $l=1$, $\eta_1>0$, and $\|\tCov\|^{-2}\|\tCov^2\| $ does not converge to zero, then
$\beta(\eta_1)$, which in this case equals $\beta_1$, cannot be estimated consistently with a standard \PLS{} algorithm.

\subsection{Examples of dependence structures}
\label{subsec:depend_ex}
In all previous theorems the ratio $\|\tCov^2\|\|\tCov\|^{-2}$
plays a crucial role. In this section some special covariance 
matrices $\tCov^2$ are studied in order to understand its behaviour.
Stationary processes considered in this section are assumed to have expectation zero and to decay exponentially, i.e., for for $c, \rho >0$ and $\gamma(0)>0$,
	\begin{align}	
		\label{eq:acfBound}
		|\gamma(t)| &\leq \gamma(0) c \exp(-\rho t), ~~ t\in \N,
	\end{align}
with $\gamma:\Z \rightarrow \R$ being the autocovariance function of the process.

Subsequently, $f(n)
\sim g(n)$ denotes  $c_1 \leq f(n)/g(n) \leq c_2$, for $n$ large,
$0<c_1<c_2$ and $f,g:\N\rightarrow \R$.
\begin{theorem}
\label{th:stationary}
Let $\left[\tCov^2\right]_{t,s}=\gamma(|t-s|)$ ($t,s=1,\ldots,n$) be the
covariance matrix of a stationary process, such that the autocovariance function $\gamma:\Z \rightarrow \R$
        satisfies (\ref{eq:acfBound}).
	Then $\|\tCov^2\| \sim n^{1/2}$ and $\|\tCov\|^2 \sim n$.
\end{theorem}

Hence, if $\tCov^2$ in model (\ref{eq:genModel}) is a covariance matrix of
a stationary process, then ignoring dependence of
observations in the \PLS{} algorithm does not affect the rate of convergence
of \PLS{} estimators, but might affect the constants. Examples of
processes with exponentially decaying autocovariances are stationary \ARMA es.

As an example of a non-stationary process we consider first-order
integrated processes. If $\{X_t\}_{t\in\Z}$ is stationary with autocovariance
function $\gamma$ satisfying (\ref{eq:acfBound}), then
$\sum_{i=1}^tX_i$ is an integrated process of order $1$. 
\begin{theorem}
\label{th:Autocovariance}
Let $\{X_t\}_{t\in\Z}$ be a stationary process with autocovariance
function $\gamma$ satisfying (\ref{eq:acfBound}). If  $\gamma(t)<0$ for some $t$, we assume additionally $\rho > \log(2c + 1)$. 
	Let $\tCov^2$ be the covariance matrix of $\sum_{i=1}^tX_i$. 
	Then $\|\tCov\|^2 \sim n^2$ and $\|\tCov^2\| \sim
        n^2$.
\end{theorem}

The lower bound on $\rho$ for negative $\gamma(t)$ ensures that no element on the
  diagonal of $\tCov^2$ becomes negative, so that $\tCov^2$ is a valid
  covariance matrix.
  
This theorem implies that the ratio $\|\tCov\|^{-2}\|\tCov^2\|$ does not
converge to zero for certain integrated processes. In particular,
combining this result with Theorems \ref{th:FirstComponentError} and \ref{th:Inconsistency} shows that the elements of the sample Krylov space $A$ and $b$, as well as $\widehat{\beta}_1$, are inconsistent, if the
dependence structure of the data can be described by an
integrated process satisfying the conditions of Theorem
\ref{th:Autocovariance}, e.g., an integrated \ARMA\ of order $(1,1,1)$. 


\section{Practical issues}
\label{sec:Estimation}
\subsection{Corrected \PLS{} estimator}
So far we considered the standard \PLS{} algorithm, showing that if
certain dependences in the data are ignored, estimation is inconsistent. Hence, it is crucial to take into account the dependence structure of
the data in the \PLS{} estimators.

Let us define $b\left(S\right) = n^{-1} X^\T S^{-2}
y$ and $A\left(S\right) = n^{-1} X^\T S^{-2} X$ for an invertible matrix $S \in \R^{n \times n}$. Furthermore, let $k_i(S)=A(S)^{i-1} b(S)$, $K_i(S)=\left[k_1(S),\ldots,k_i(S)\right] \in \R^{k \times i}$ and $\widehat{\beta}_i(S) = K_i(S)\left\{ K_i(S)^\T A(S) K_i(S) \right\}^{-1} K_i(S)^\T b(S)$ ($i=1,\ldots,k$).

For $S = I_n$ this yields a standard \PLS{} estimator. If
$S = \tCov$, the temporal dependence matrix, then $b\left(\tCov\right)$ and
$A\left(\tCov\right)$ are square-root consistent estimators of $Pq$ and
$\Sigma^2$, respectively, with the mean squared error independent of
$\tCov$, which follows from Theorem \ref{th:FirstComponentError}. Hence, the resulting $\widehat{\beta}_i(\tCov)$ is also a
consistent estimator of $\beta_i$ and Theorem \ref{th:ExpBlanchard} shows that $\beta(\eta_1)$ can be estimated consistently by early stopping as well. This
procedure is equivalent to running the \PLS{} algorithm on
$\tCov^{-1}y$ and $\tCov^{-1}X$, that is, with the temporal
dependence removed from the data.

In practice the true covariance matrix $\tCov^2$ is typically unknown
and is replaced by a consistent estimator $\widehat{\tCov}^2$. We call the estimator
$\widehat{\beta}_i(\widehat{\tCov})$ the corrected \PLS{} estimator. 
The next theorem shows that, given a consistent estimator of $\tCov^2$, the population Krylov space and $\beta(\eta_1)$ can be estimated consistently.
\begin{theorem}
\label{th:consEstimator}
Let $\widehat{\tCov}^2$ be an estimator for $\tCov^2$ that is invertible for all $n \in \N$ and
$\left\|\tCov \widehat{\tCov}^{-2} \tCov  - I_n \right\|_\calL = O_p(r_n)$, where $r_n$ is some sequence of positive numbers such that $r_n\rightarrow0$ as $n\rightarrow\infty$. Then 
\[
	\|A(\widehat{\tCov}) - \Sigma^2\|_\mathcal{L} = O_p(r_n), \quad
	\|b(\widehat{\tCov})-P q\| = O_p(r_n).
\]
Moreover, with probability at least $1-\nu$, $\nu\in(0,1]$
$$\|\widehat{\beta}_{a^\ast}(\widehat{\tCov}) - \beta(\eta_1)\| = O(r_n),$$ where the definition of $a^*$ in (\ref{eq:discrepancy0}) is updated by replacing $A$, $b$ and $\widehat{\beta}_i$ by $A(\widehat{\tCov})$, $b(\widehat{\tCov})$ and $\widehat{\beta}_i(\widehat{\tCov})$, respectively.
\end{theorem}

Theorem \ref{th:consEstimator} states that if a consistent estimator of the covariance matrix $V^2$ is available, then the elements of the population Krylov space $A$, $b$, as well as the coefficient $\beta(\eta_1)$, can be consistently estimated by  $A(\widehat{\tCov})$, $b(\widehat{\tCov})$ and $\widehat{\beta}_{a^\ast}(\widehat{\tCov})$. The convergence rate of these estimators is not faster than that of $\widehat{\tCov}^2$. For example, if 
the temporal dependence in the data follows some parametric model, then parametric rates of $n^{-1/2}$ are also achieved for $A(\widehat{\tCov})$, $b(\widehat{\tCov})$ and $\widehat{\beta}_{a^\ast}(\widehat{\tCov})$. Estimation of $\tCov^2$ by some nonparametric methods, e.g., with a banding or tapering approach, leads to a slower convergence rates: see \citet{Bic08} or \citet{Wu12}. Similar results are well-known in the context of linear regression. For example, Theorem 5.7.1 in \citet{Fuller1997} shows that the convergence rate of feasible generalized least squares estimators is the same as that of the estimator for the covariance matrix of the regression error.

%


\subsection{Estimation of covariance matrices}
\label{sec:cov.matrix.estimation}
To obtain the corrected \PLS{} estimator, some consistent estimator of
$\tCov^2$ based on a single realisation of the process is
necessary. In model
(\ref{eq:genModel}) the dependence structure over the observations of
$X$ is the same as that of $y$ and $\tCov$ can be estimated from $y$ alone.

If $\tCov^2$ is the autocovariance matrix of a stationary process, it can be estimated both parametrically and nonparametrically. 
Many stationary processes can be sufficiently well approximated by an autoregressive moving average process,  see
\citet{bBrock}, Chapter 4.4. Parameters of autoregressive moving average processes are estimated either by Yule--Walker or maximum likelihood estimators, both attaining parametric rates. 
Another approach is to band or taper the empirical
autocovariance function of $y$ \citep{Bic08,Wei09,Wu12}. These nonparametric estimators are very flexible, but are computationally intensive and have slower convergence rates. 

If $y$ is an integrated processes of
order one, then $\tCov^2$ can easily be derived from the
covariance matrix estimator of the corresponding stationary process.

\section{Simulations}
\label{sec:Simulation}
To verify small sample performance of the \PLS{} algorithm under
dependence we consider the following simulation setting. To illustrate consistency we
choose three sample sizes $n \in\{ 250,500,2000\}$. In the latent
variable model (\ref{eq:genModel}) we set $k=20$, $l = 1,5$ and take the elements of $P$ to be  independent identically distributed Bernoulli random
variables with success probability $0.5$. 
Elements of the vector $q$ are $q_i=5 \, i^{-1}$ ($i=1,\ldots,l$), in order to control the importance of the
different latent variables for $y$.
The random variables $N_{1,1}$, $F_{1,1}$ and $f_1$ are
taken to be standard normally distributed.
The parameter $\eta_2$ is chosen to get the signal to noise ratio in $y$ to be two and $\eta_1$ is set so that the signal to noise ratio in $X$ is $0.5$.
Three matrices $\tCov^2$ are considered: the identity
matrix, the covariance matrix of an autoregressive process of the first order with coefficient $0.9$ and the covariance matrix of an autoregressive integrated moving average process of order $(1,1,1)$ with both parameters set to $0.9$. 

First, we ran the standard \PLS{}
algorithm on the data with the three aforementioned dependence structures to highlight the effect of the ignored dependence in the data. Next, we studied the performance of our corrected \PLS{} algorithm applied to non-stationary data. Thereby, the covariance matrix of the autoregressive moving average process has been estimated parametrically, as discussed in Section \ref{sec:cov.matrix.estimation}. A nonparametric estimation of this covariance matrix has lead to qualitative similar results.



\begin{figure}[t!]
\captionsetup{width=0.9\textwidth}
\begin{minipage}{.5\textwidth}
  	\centering
	\figurebox{16pc}{16pc}{}[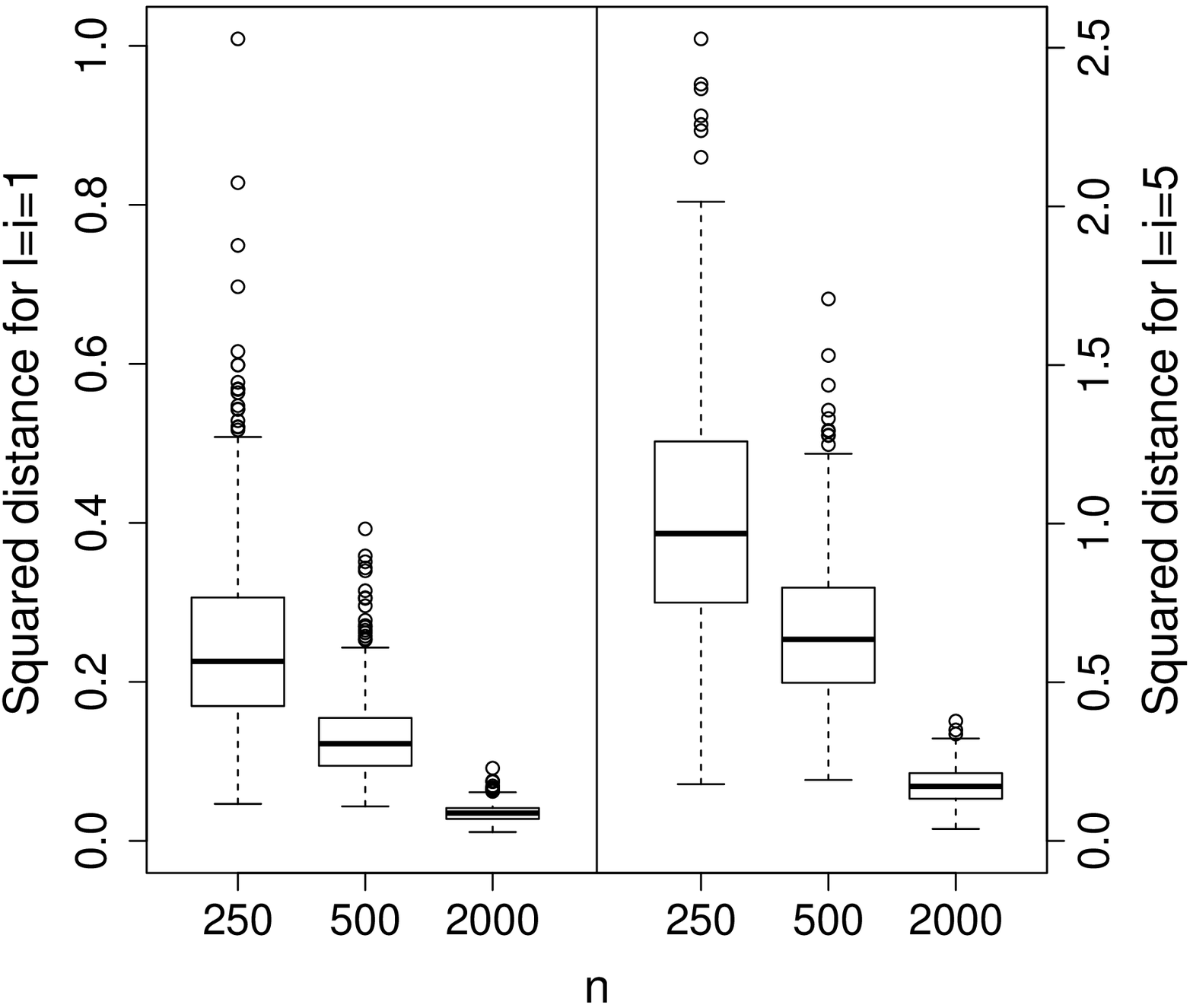]
\end{minipage}%
\begin{minipage}{.5\textwidth}
  	\centering
	\figurebox{16pc}{16pc}{}[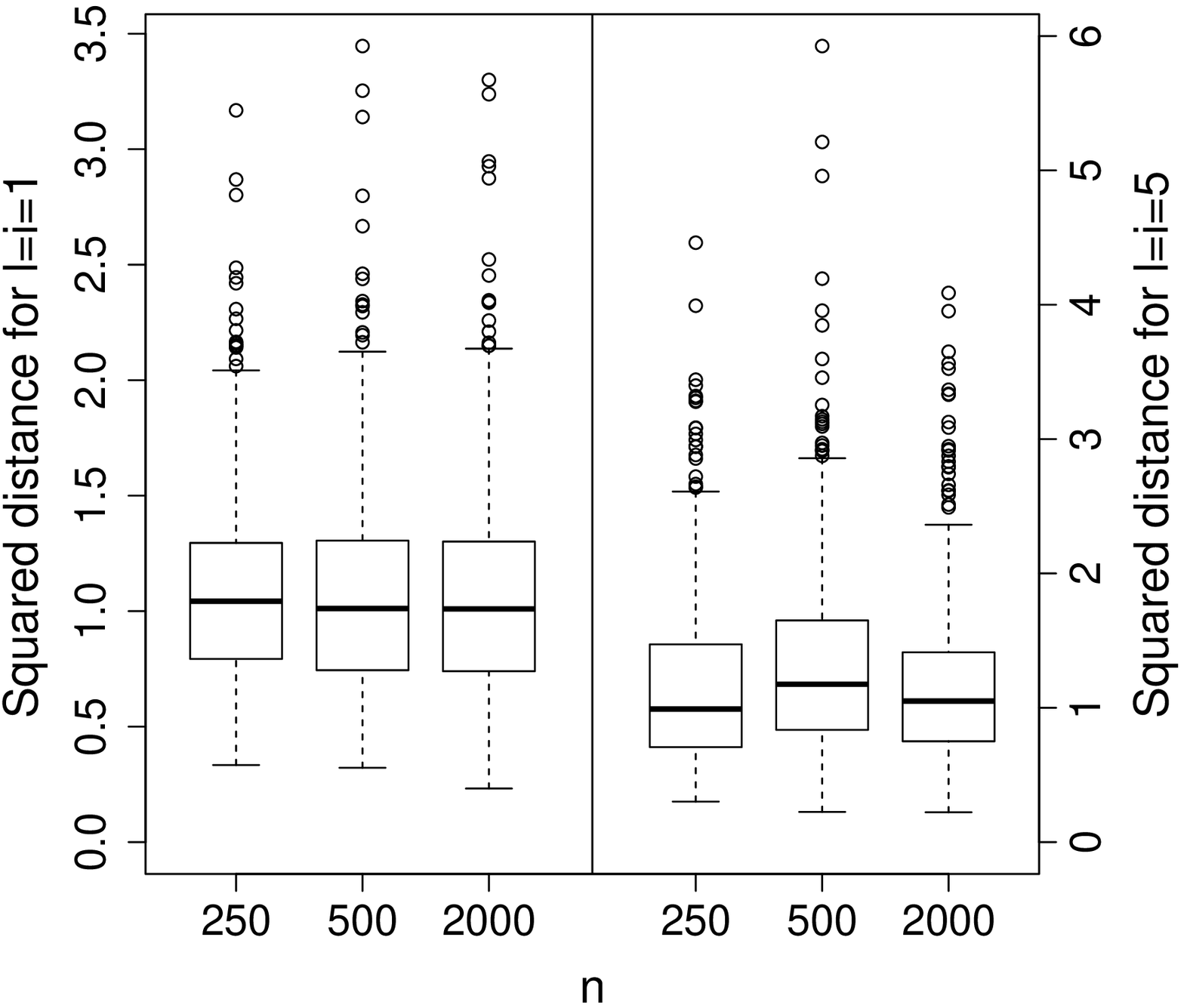]
\end{minipage}
\begin{minipage}{.5\textwidth}
  	\centering
	\figurebox{16pc}{16pc}{}[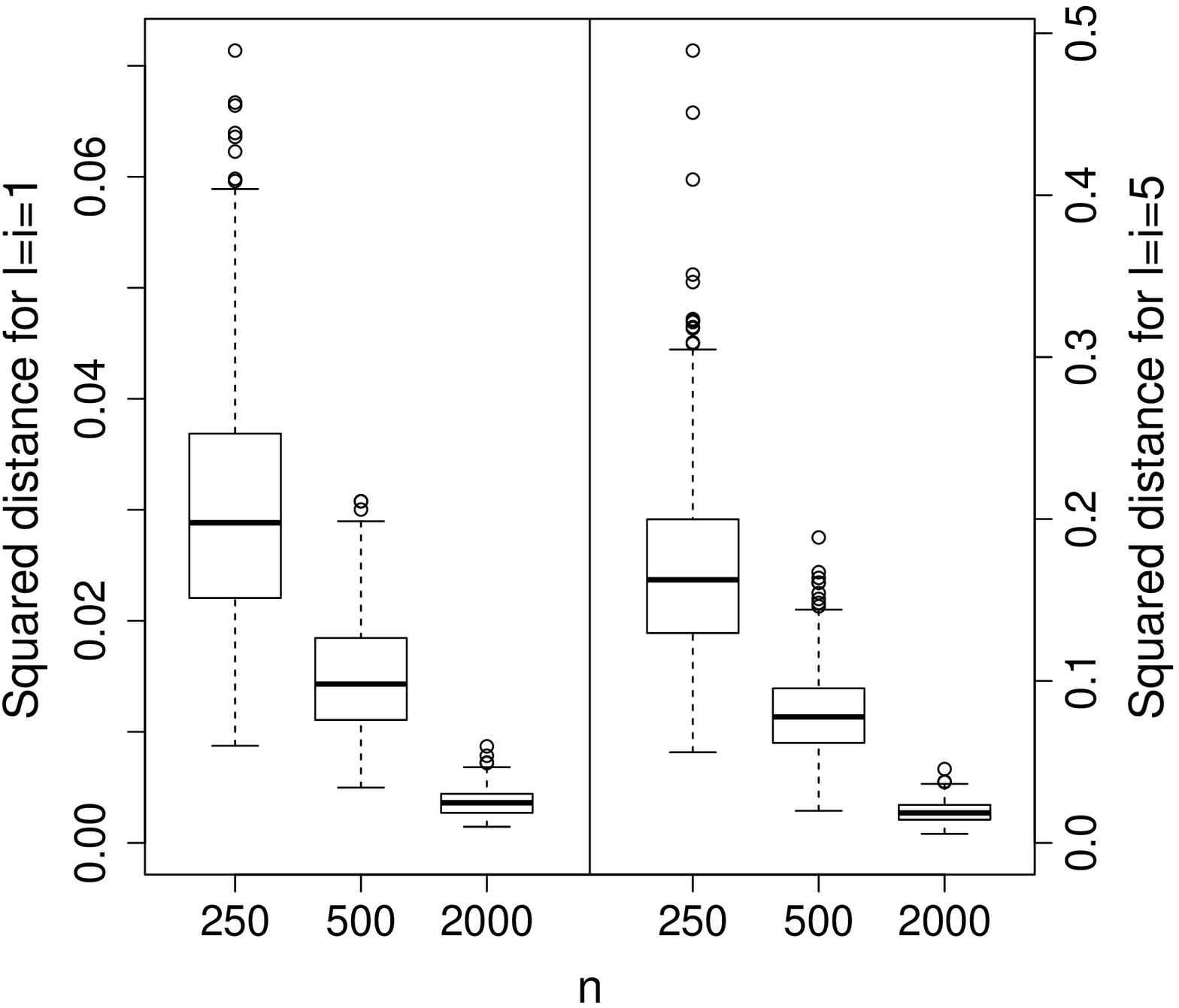]
\end{minipage}
\begin{minipage}{.5\textwidth}
  	\centering
	\figurebox{16pc}{16pc}{}[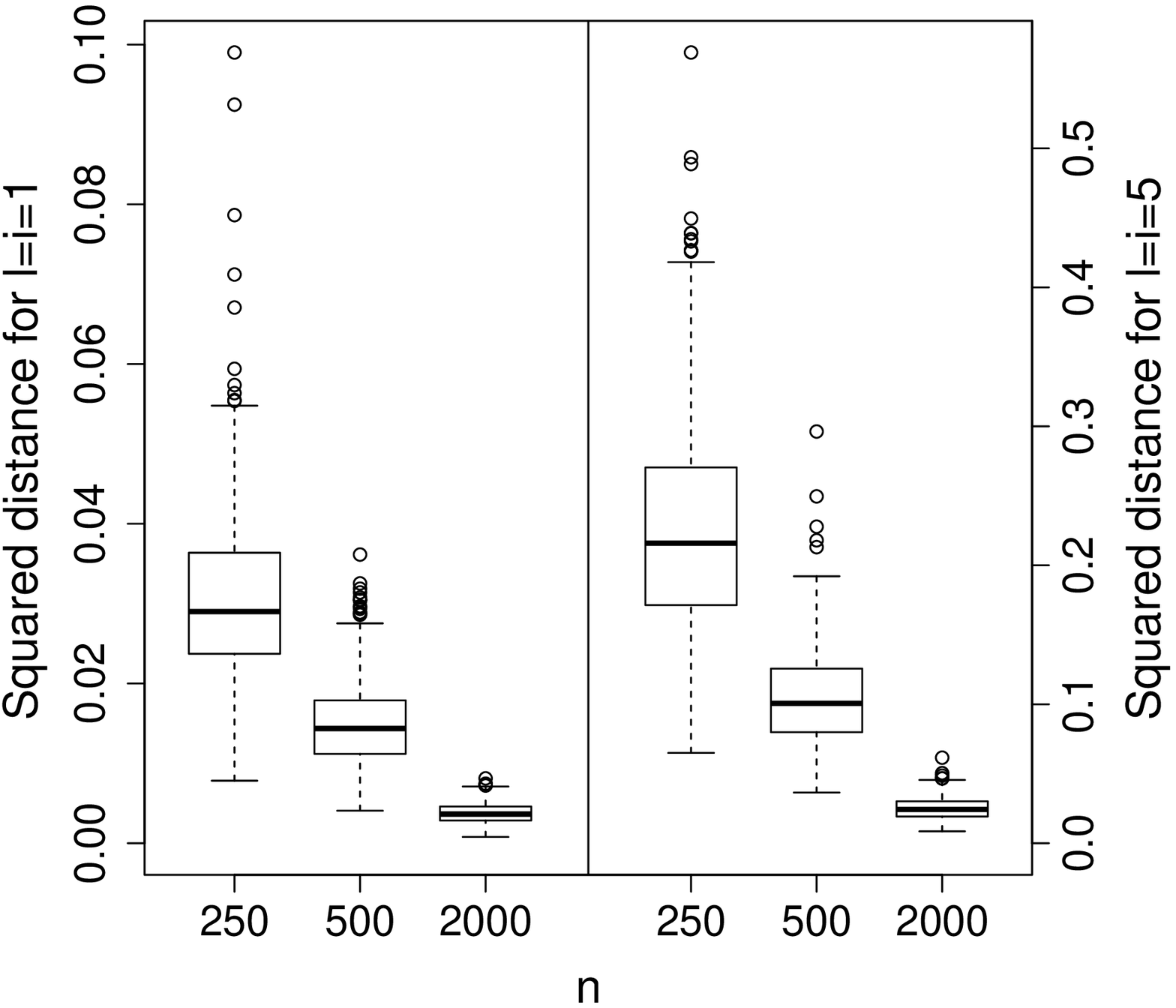]
\end{minipage}
\caption{Squared distance of partial least squares estimators $\widehat{\beta}_i$ and $\beta(\eta_1)$ in 500 Monte Carlo samples. First three boxplots in each panel correspond to $l=i=1$, the latter three to $l=i=5$. The dependence structures are: first order autoregressive (top left), autoregressive integrated moving average of order (1,1,1) (right) and independent, identically distributed (bottom left). The standard partial least squares (top and bottom left) and corrected partial least squares (bottom right) have been employed.\label{fig:EMSE}}
\end{figure}

The boxplots in Figure \ref{fig:EMSE} show the squared distance of $\widehat{\beta}_i$ and $\beta(\eta_1)$ in $500$ Monte Carlo replications. Two
cases are shown in one panel: the model has just one latent component and $\widehat\beta_1$ is
considered, i.e., $l=i=1$ and the model has five latent components and the
squared distance of $\widehat\beta_5$ to $\beta(\eta_1)$ is studied, i.e., $l=i=5$.  

We observe that the mean squared error
of $\widehat\beta_i$ obtained with the standard partial least squares converges to zero for
autoregressive and independent data with the growing sample size. However, an autoregressive
dependence in the data leads to a somewhat higher mean squared error, compare the top and bottom left panels. If the data follow an autoregressive integrated moving average process and this is ignored in the \PLS{} algorithm,
then the mean squared error of $\widehat\beta_i$ converges to some
positive constant, see the top right boxplots. Taking into account these non-stationary dependencies in the corrected partial least squares leads to consistent estimation, similar to the independent data case, compare the bottom left and right panels.


We conclude
that if the observations are dependent, corrected \PLS{} improves
estimation: in case of stationary dependence the mean squared
error is reduced and in case of non-stationary dependence the estimation
becomes consistent.

\section{Application to Protein Dynamics}
\label{sec:Data}
Proteins
\begin{figure}[b!]
\captionsetup{width=0.9\textwidth}
\begin{minipage}{.5\textwidth}
  	\centering
	\figurebox{16pc}{16pc}{}[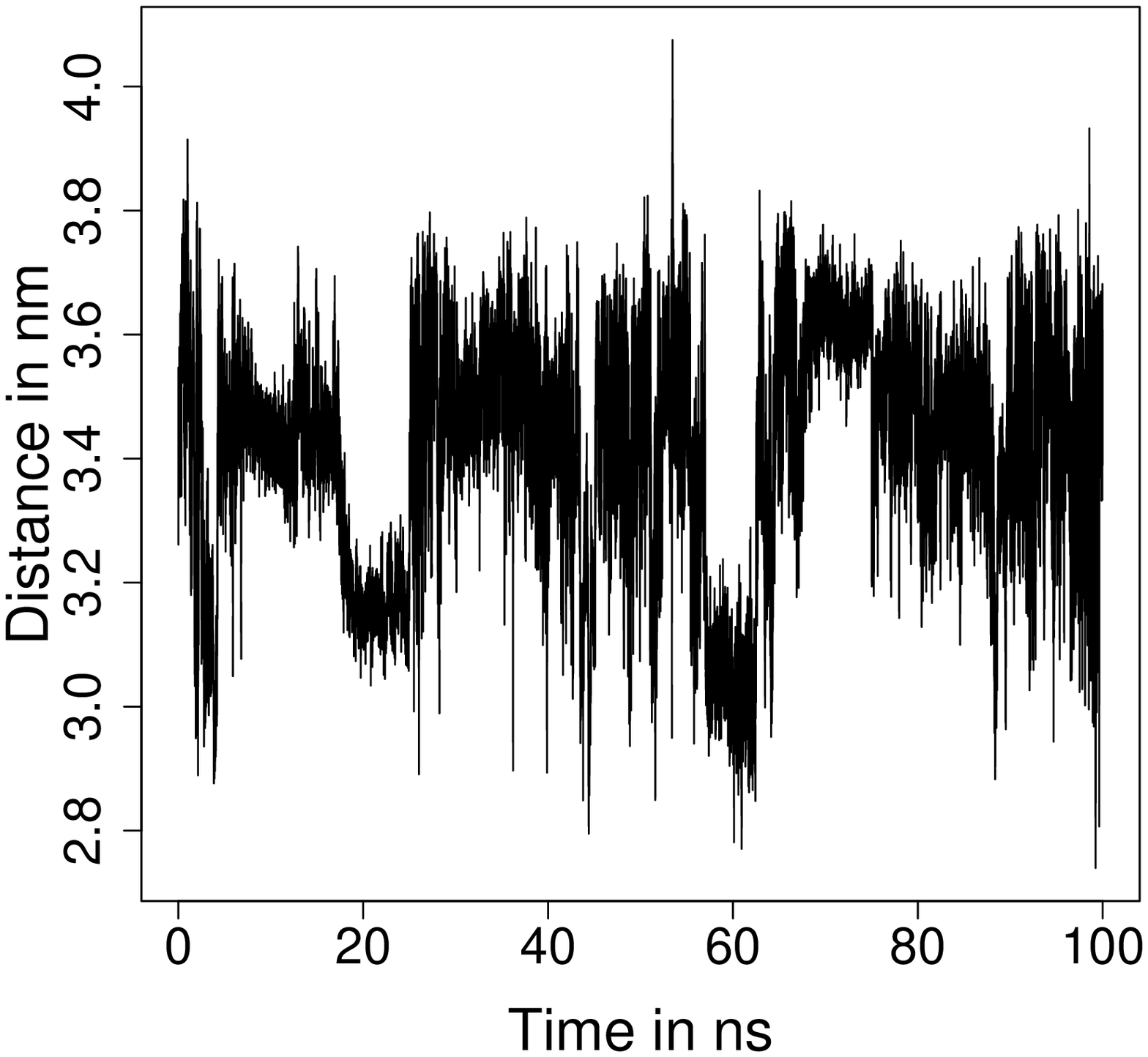]
\end{minipage}
\begin{minipage}{.5\textwidth}
  	\centering
	\figurebox{16pc}{16pc}{}[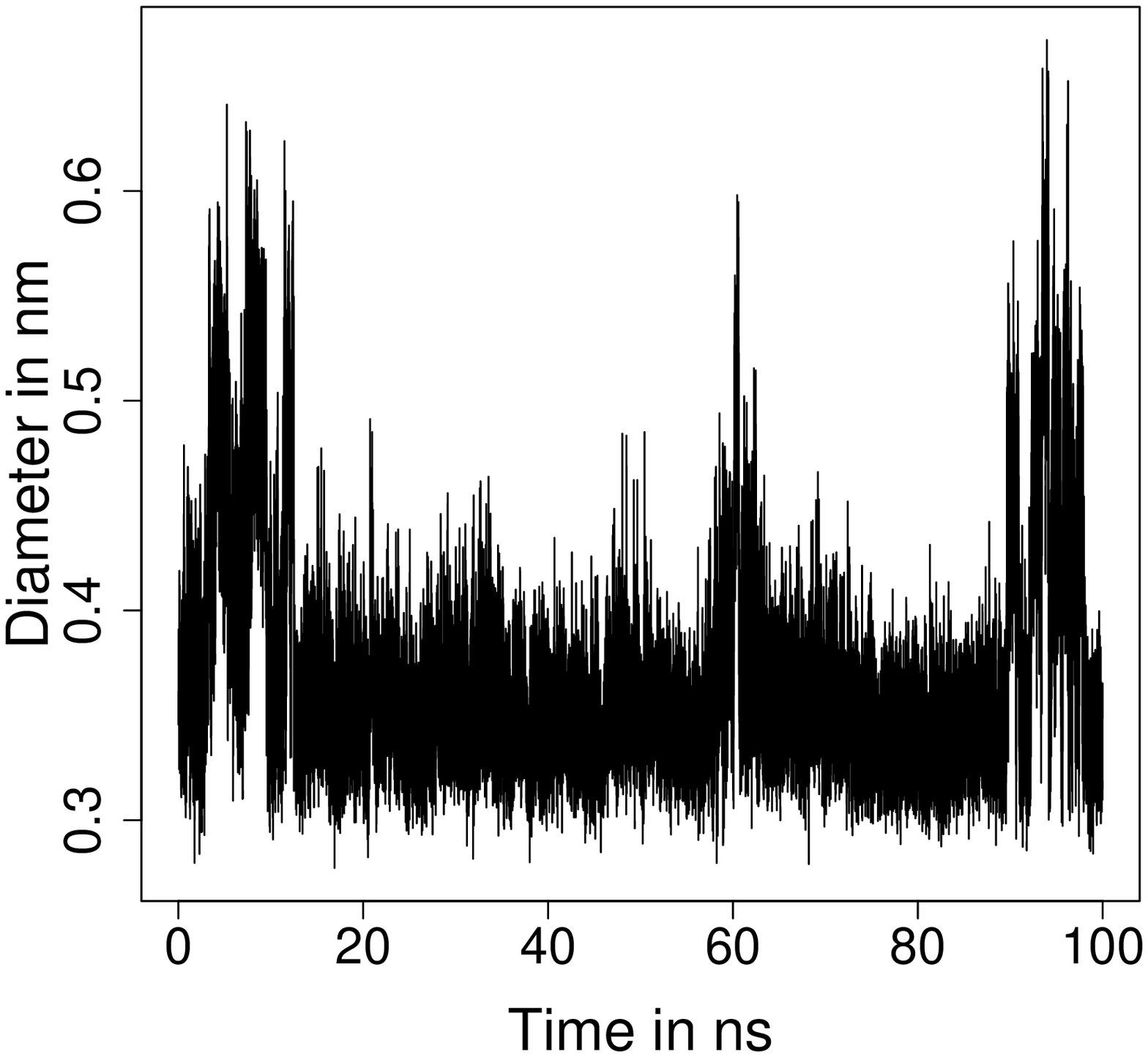]
\end{minipage}
\caption{Distance between the first backbone atom and the first centre of mass of aquaporine (left) and the opening diameter over time (right). \label{fig:AQY1data}}
\end{figure} fulfil their biological function through particular movements, see \citet{Henz07}, so a key step in understanding protein functions
is a detailed knowledge of the underlying dynamics. Molecular dynamics simulations  \citep{BdG1998} are routinely used to study
the dynamics of biomolecular systems at atomic detail on timescales of
nanoseconds to microseconds. 
Although in principle allowing
to directly address function-dynamics relationships, analysis is
frequently hampered by the large dimensionality of the protein
configuration space, rendering it non-trivial to identify collective
modes of motion that are directly related to a functional property of
interest.

\citet{KrivGroot12} have shown that \PLS{} helps
to identify a hidden relation between
atom coordinates of a protein and a functional parameter of
interest, yielding robust and parsimonious
solutions, superior to principal component regression.
In this work we look at a protein
studied in the aforementioned paper: the water channel aquaporine as found in the yeast Pichia
pastoris. This is a gated channel, i.e., the diameter of the opening
can change, controlling the flow of water into the cell. 
We aim to study which collective
motions of protein atoms influence the diameter $y_t$ of the channel at time
$t$, measured by the distance of two centres of mass of the
residues of the protein which characterize the opening. For the
description of the protein dynamics we use an inner model, i.e. at
each point in time we calculate the Euclidean distance $d$ of each backbone atom of
the protein and a set of certain four fixed base points. 
We denote the $p=739$ atoms by $A_{t,1},\ldots,A_{t,p} \in \R^3$, the fixed base points by $B_1,\ldots,B_4 \in \R^3$ and take
\[
	X_t = \left\{
		d(A_{t,1},B_1),\ldots,d(A_{t,p},B_1),d(A_{t,1},B_2),\ldots,d(A_{t,p},B_4)\right\}^\T\in\R^{4p}.
\]
The available timeframe has a length of $100$ ns split into
$n=20\,000$ equidistant points of observation. \citet{KrivGroot12}
found that a linear relationship between $X$ and $y$ can be
assumed. Additionally, these data seem not to contradict model (\ref{eq:genModel}).
Taking a closer look at the data reveals that both $y_t$ and
$X_{t,i}$ ($i=1,\ldots,4p$) are non-stationary time series, see Figure \ref{fig:AQY1data}.
For the calculation of $\widehat{\tCov}^2$ we used the banding approach
  mentioned in Section \ref{sec:cov.matrix.estimation} and found the results to
  be very
  similar to a simple autoregressive integrated moving average process with parameters (3,1,1) and corresponding coefficients
  $(0.1094,0.0612,0.0367, -0.9159)$. Autoregressive integrated moving average models have been employed before to study protein time series \citep{Alak}.

\begin{figure}[t!]
\captionsetup{width=0.9\textwidth}
\begin{minipage}{.5\textwidth}
  	\centering
	\figurebox{16pc}{16pc}{}[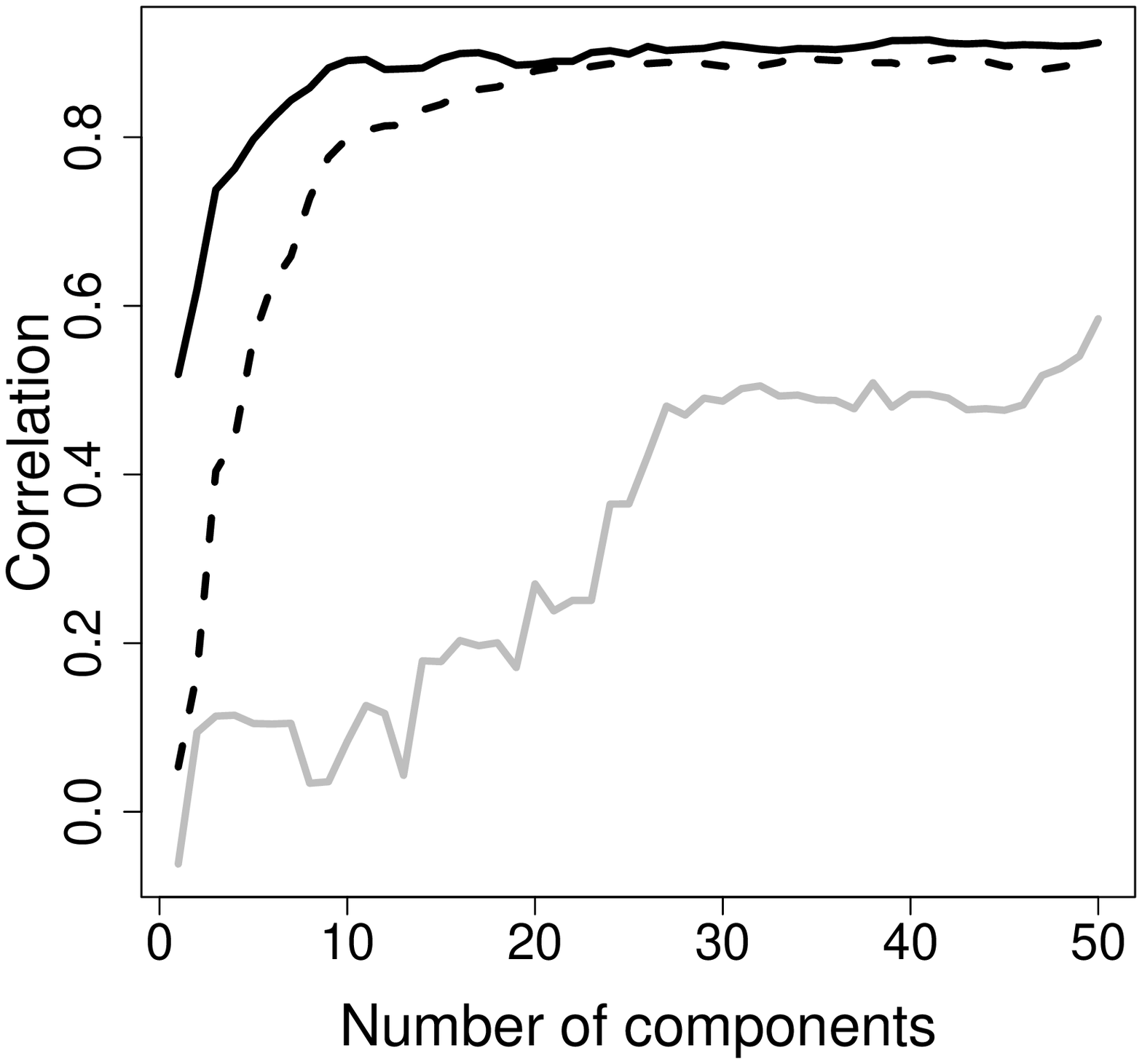]
\end{minipage}
\begin{minipage}{.5\textwidth}
  	\centering
	\figurebox{16pc}{16pc}{}[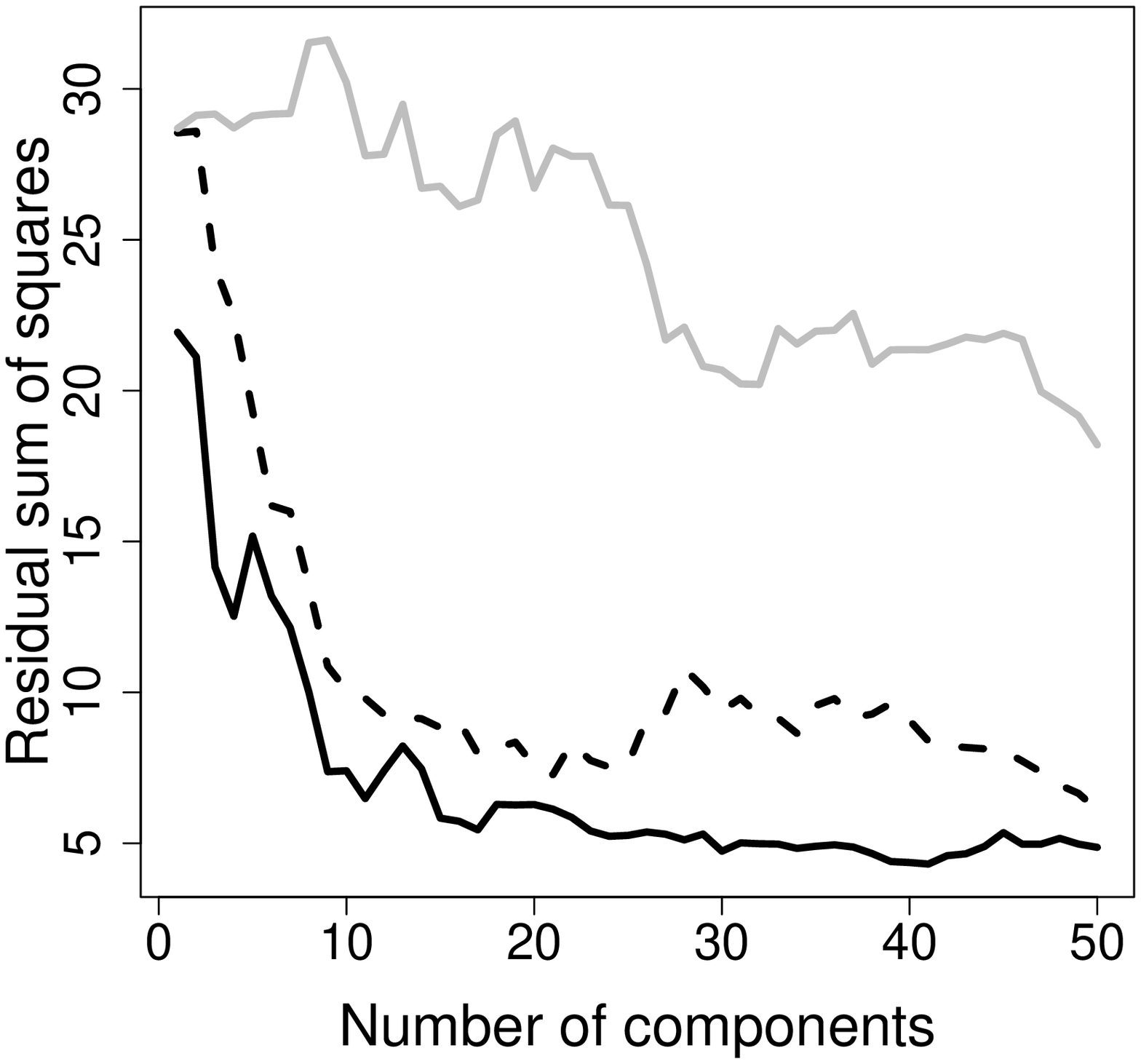]
\end{minipage}
\caption{Correlation (left) and residual sum of squares (right) of the predicted opening diameter and the real data on the test set. Compared methods are principal component regression (grey), corrected partial least squares (black, solid) and partial least squares (black, dashed). \label{fig:AQY1}}
\end{figure}

To validate our estimators, we used the following procedure.
 First, the data were split into two equal parts and the models were build
 on the first half. Then the prediction was done on the
 test set consisting of the second half of the data and was compared to $y_t$ from the test set. To measure the
 accuracy of the prediction we used the Pearson correlation
 coefficient common in the biophysics community and the
 residual sum of squares, both shown in Figure \ref{fig:AQY1}. The \PLS{}
 estimator clearly outperforms principal components regression.
  The corrected \PLS{}
 algorithm, which takes temporal dependence into account, delivers better prediction than standard \PLS{}. The
 improvement is strongly present in the first components.

High predictive power of the first corrected \PLS{} components is particularly
        relevant for the interpretation of the underlying protein dynamics. \citet{KrivGroot12} established that the first \PLS{}
        regression coefficient $\widehat{\beta}_1$ corresponds to the
        so-called ensemble-weighted maximally correlated mode of motion contributing most to the fluctuation in the response $y$.
Altogether, due to the low dimensionality, corrected \PLS{} greatly facilitates the
interpretation of the underlying relevant dynamics, compared to \PLS{}
and principal component regression, where many more components are required to obtain the same predictive power.

\section*{Acknowledgements}
The authors are grateful to the reviewers for their valuable comments which helped to improve the article. The support of the German Research Foundation is gratefully acknowledged.

\section*{Supplementary Material}
Supplementary Material available at {\it Biometrika} online includes all technical details. 

\bibliography{bibliography}
\biblist
\end{document}


\nolinenumbers


\markboth{M. Singer, T. Krivobokova, A. Munk \and B. de Groot}{Partial least squares for dependent data}

\title{Supplementary material for \\ 
Partial least squares for dependent data}

\author{MARCO SINGER, TATYANA KRIVOBOKOVA, AXEL MUNK}
\affil{Institute for Mathematical Stochastics, Goalgorithm2e.styldschmidtstr. 7, 37077 G\"ottingen, Germany
\email{msinger@gwdg.de} 
\email{tkrivob@uni-goettingen.de}
\email{munk@math.uni-goettingen.de}}
\author{\and BERT DE GROOT}
\affil{Max Planck Institute for Biophysical Chemistry, Am Fassberg 11, 37077 G\"ottingen, Germany 
\email{bgroot@gwdg.de}}

\maketitle

\section{Proofs}
\label{sec:proofs}
\subsection{Derivation of the population \PLS{} components}
Let denote $K_i \in \R^{k \times i}$ the matrix representation of a base for $\mathcal{K}_i(\Sigma^2,\spCC)$ . Then
\begin{align*}
	\sum\limits_{t=1}^n \E\left(y_t - X^\T_t K_i \alpha\right)^2= \sum\limits_{t=1}^n [\tCov^2]_{t,t} \left(
		\|q\|^2 + \eta_2^2 - 2 \alpha^\T K_i^\T \spCC + 
		\alpha^\T K_i^\T \Sigma^2 K_i \alpha
	\right).
\end{align*}
Minimizing this expression with respect to $\alpha \in \R^i$ gives
$
	K_i^\T \Sigma^2 K_i \alpha = K_i \spCC.
$ 
Since the matrix $K_i^\T \Sigma^2 K_i$ is invertible, 
we get the least squares fit $\beta_i$ in Section 2. 

Assume now that the first $i<a$ \PLS{} base vectors $w_1,\ldots,w_i$ have been calculated and consider for $\lambda\in \R$ the Lagrange function
\begin{align*}
	\sum\limits_{t,s=1}^n \Cov\left(y_t - X_t^\T \beta_i, X_s^\T w\right) -  \lambda(\|w\|^2 -1)=
	w^\T\left(
		\spCC - \Sigma^2\beta_{i}
	\right)\sum\limits_{t,s=1}^n [\tCov^2]_{t,s}  - \lambda (\|w\|^2-1).
\end{align*}
Maximizing with respect to $w$ yields
$$
	w_{i+1} = (2\lambda)^{-1}\left(
		\spCC - \Sigma^2 \beta_i 
	\right)
	\sum\limits_{t,s=1}^n [\tCov^2]_{t,s} 
	\propto \spCC - \Sigma^2 \beta_i.
$$
Since $\beta_i \in \mathcal{K}_i(\Sigma^2,\spCC)$, we get $w_{i+1} \in \mathcal{K}_{i+1}(\Sigma^2,\spCC)$ and $w_{i+1}$ is orthogonal to $w_1,\ldots,w_i$.

\subsection{Proof of Theorem 1}
First consider
\begin{align*}
	\notag
	\E\left(\left\| b - \spCC
	\right\| ^2\right) = &\E\left[\left\|
		\frac{1}{\|\tCov\|^2}\left\{
			(P N^\T + \eta_1 F^\T)\tCov^2 Nq + \eta_2(P N^\T  + \eta_1F^\T)\tCov^2f
		\right\} - \spCC
	\right\|^2\right]\\
	=& \left\{\E\left( \left\|\frac{1}{\|\tCov\|^2}P N^\T \tCov^2 N q - Pq \right\|^2\right) + 
	\frac{\eta_2^2}{\|\tCov\|^4}\E\left(\left\| P N^\T \tCov^2 f \right\|^2\right)\right\}\\ 	
	&+ \frac{\eta_1^2}{\|\tCov\|^4}\left\{ \E\left(
	\left\| 
		F^\T \tCov^2 N q 
	\right\|^2\right) + \eta_2^2\E\left( \left\| F^\T \tCov^2 f \right\|^2\right) \right\} = S_1+S_2,
\end{align*}
due to the independence of $N$, $F$ and $f$. It is easy to see that
\begin{equation*}
S_2= \frac{\|\tCov^2\|^2}{\|\tCov\|^4} \eta_1^2 k \left(
		 \|q\|^2 +  \eta_2^2
	\right).
\end{equation*}
Furthermore, with notation $A_0=N^\T \tCov^2 N$ we get
\begin{equation*}
	S_1 = \frac{1}{\|\tCov\|^4}\E \left(q^\T 
		A_0 P^\T P A_0
	q\right)
	- \left\|P q \right\|^2 + \frac{\eta_2^2}{\|\tCov\|^4}\E\left(\left\| P N^\T \tCov^2 f \right\|^2\right).
\end{equation*}
Consider now $\E \left(q^\T 
		A_0 P^\T P A_0
	q\right)$ as a quadratic form with respect to the matrix
      $P^\T P$. {Denote $\kappa=\E\left(N_{1,1}^4\right) - 3$}. First, 
$
	\E\left(A_0 q\right) = \E\left(N^\T \tCov^2 N q\right) = \|\tCov\|^2 q
$
and
\begin{align*}
	\Var(A_0q) 
	&= \left[
			\sum\limits_{a,b=1}^l q_a q_b \sum
			\limits_{t,s,u,v=1}^n  
			\tCov_u^\T \tCov_s
			\tCov_t^\T \tCov_v \E(N_{s,i} N_{u,a}N_{t,j}N_{v,b})
		\right]_{i,j=1}^l - \|\tCov\|^4 q q^\T\\
	& = \left[
		q_iq_j\|\tCov\|^4 + \left(q_iq_j + \delta_{i,j}\|q\|^2\right) \|\tCov^2\|^2 + \kappa \sum\limits_{t=1}^n \|\tCov_t\|^4 \delta_{i,j} 
		q_i^2\right]_{i,j=1}^l - \|\tCov\|^4 q q^\T \\
	& =  \|\tCov^2\|^2 \left(qq^\T + \|q\|^2I_l\right) + \kappa \sum\limits_{t=1}^n \|\tCov_t\|^4 \mathrm{diag}\left( q_1^2,\ldots,q_l^2\right),
\end{align*}
{where $\mathrm{diag}(v_1,\dots,v_l)$ denotes the diagonal matrix with entries $v_1,\dots,v_l \in \R$ on its diagonal and $\delta$ is the Kronecker delta}.
In the second equation we made use of
$	\E\left( N_{s,i}N_{u,a}N_{t,j}N_{v,b}\right) = \delta_{i,a}\delta_{j,b}\delta_{s,u}\delta_{t,v}+
		\delta_{i,b}\delta_{j,a}\delta_{s,v}\delta_{t,u}+	
		\delta_{i,j}\delta_{a,b}\delta_{t,s}\delta_{u,v} +
		\kappa\, \delta_{t,s}\delta_{s,u}\delta_{u,v}\delta_{i,j}\delta_{j,a}\delta_{a,b}
$.
Hence,
\begin{align*}
	\frac{1}{\|\tCov\|^4} \E \left(
		q^\T A_0 P^\T P A_0
	q \right)
	 =& \frac{1}{\|\tCov\|^4} \trace\left\{
		P^\T P \Var\left(A_0q\right)\right\}
	 - \frac{1}{\|\tCov\|^4} \E\left(q^\T A_0 \right) P^\T P
	\E\left(A_0 q\right)\\
	=& \frac{\|\tCov^2\|^2}{\|\tCov\|^4} \left(
		q^\T P^\T P q + \|P\|^2 \|q\|^2
	\right) + q^\T P^\T P q+\kappa\sum\limits_{t=1}^n \frac{\|\tCov_t\|^4}{\|\tCov\|^4} 
		\sum\limits_{i=1}^l \|P_i\|^2 q_i^2.
\end{align*}
The remaining term in $S_1$ follows trivially, proving the
result. $\E\|\Sigma^2-A\|^2$ is obtained using similar calculations.
\eop

\subsection{Proof of Theorem 2}
\begin{lemma}
\label{lem:sqrRootConcentration}
Assume that for $\nu \in (0,1]$ and some constants $\delta,\epsilon>0$ it holds that $\Prob\left(\|A - \Sigma^2\|_\calL \leq \delta\right) \geq 1- \nu/2$ and $\Prob\left(\|b - P q \| \leq \epsilon\right) \geq 1- \nu/2$. Then each of the inequalities
	\begin{align*}
		\|A^{1/2} - \Sigma\| &\leq  2^{-1} \delta \|\Sigma^{-1}\|  \{1 + o(1)\},\\ 
		\|A^{-1/2} b - \Sigma^{-1} Pq\| &\leq  \epsilon \|\Sigma^{-1}\|_\calL  + 2^{-1} \delta(\|Pq\| + \epsilon)
		 \|\Sigma^{-2}\|\|\Sigma^{-1}\|
	\left\{ 1 + o(1)\right\}
	\end{align*}
	hold with probability at least $1-\nu/2$.
\end{lemma}
{\it{Proof}}: We show the result by using the Fr\'echet-derivative for functions $F:\R^{k \times k} \rightarrow \R^{k \times k}$. Due to the fact that $\eta_1>0$ it holds that $\Sigma^2$ is positive definite and thus invertible.

It holds due to \citet{Hig2008}, Problem 7.4, that $F'(\Sigma^2)B$ for an arbitrary $B \in \R^{k \times k}$ is given as the solution in $X \in \R^{k \times k}$ of $B = \Sigma X + X \Sigma$, i.e. due to the symmetry and positive definitiness of $\Sigma$ we have $F'(\Sigma^2)B = 2^{-1}\Sigma^{-1}B$. We take the  orthonormal base $\{E_{i,j}, i,j=1,\dots,k\}$ for the space $(\R^{k \times k},\|\cdot\|)$ with $E_{i,j}$ corresponding to the matrix that has zeros everywhere except at the position $(i,j)$, where it is one. 
The Hilbert-Schmidt norm $\|F'(\Sigma^2)\|_{HS}$ is 
\[
	\|F'(\Sigma^2)\|_{HS}^2 = 4^{-1}\sum_{i,j=1}^k\|\Sigma^{-1} E_{i,j}\|^2 = 4^{-1}\sum\limits_{i,j=1}^k [\Sigma^{-1}]_{i,j}^2 = 4^{-1} \|\Sigma^{-1}\|^2.
	\]
This yields with the Taylor expansion for Fr\'echet-differentiable maps
\begin{align*}
\|A^{1/2} - \Sigma\|_\mathcal{L} \leq \| F'(\Sigma)(A - \Sigma^2)\| + o(\|A - \Sigma^2\|)
\leq 
	2^{-1} \|\Sigma^{-1}\| \delta \{1 + o(1)\}.
\end{align*}
For the second inequality we see first that
\begin{equation}
\label{eq:Ainverseequatliy}
\|A^{-1/2} b - \Sigma^{-1} Pq\| \leq \epsilon \|\Sigma^{-1}\|_\calL  + \left\|
		(A^{-1/2} - \Sigma^{-1})b
	\right\|.
\end{equation}
The Fr\'echet-derivative of the map $F:\R^{k \times k} \rightarrow \R^{k \times k}, A \mapsto A^{-1/2}$ is $F'(\Sigma^2)B = -2^{-1} \Sigma^{-2} B \Sigma^{-1}$ and
\[
	\|F'(\Sigma^2)\|_{HS}^2 = 4^{-1} \sum\limits_{i,j=1}^k \|\Sigma^{-2} E_{i,j} \Sigma^{-1}\|^2 \leq 4^{-1} \|\Sigma^{-2}\|^2\|\Sigma^{-1}\|^2.
\]
Here we used the submultiplicativity of the Frobenius norm with the Hadamard product of matrices.
Thus we get via Taylor's theorem
\[
	\|A^{-1/2} - \Sigma^{-1}\| \leq 2^{-1} \|\Sigma^{-2}\|\|\Sigma^{-1}\|\|A - \Sigma^2\| + o(\delta).
\]
Plugging this into (\ref{eq:Ainverseequatliy}) yields
\[
	\|A^{-1/2} b - \Sigma^{-1} Pq\| \leq  \epsilon \|\Sigma^{-1}\|_\calL  + 2^{-1} \delta(\|Pq\| + \epsilon)
		 \|\Sigma^{-2}\|\|\Sigma^{-1}\|
	\left\{ 1 + o(1)\right\},
\]
where we used that $\|b\| \leq  \|Pq\| + \epsilon$.
\eop
\ \\
{\bf{Equivalence of conjugate gradient and partial least squares:}} We denote $\tilde{A} = A^{1/2}$ and $\tilde{b} = A^{-1/2} b$. The partial least squares optimization problem is
\[
	\min\limits_{v \in \mathcal{K}_i(A,b)} \|y - X v\|^2,
\]
whereas the conjugate gradient problem studied in \citet{Nemirovskii86} is
\begin{equation}
\label{eq:nem.cg}
		\min\limits_{v \in \mathcal{K}_i(\tilde{A}^2,\tilde{A}\tilde{b})} \|\tilde{b} - \tilde{A} v\|^2.
\end{equation}
It is easy to see that the Krylov space $\mathcal{K}_i(\tilde{A}^2,\tilde{A}\tilde{b}) = \mathcal{K}_i(A,b)$ ($i=1,\dots,k$). We have
\[
	\arg\min\limits_{v \in \mathcal{K}_i(\tilde{A}^2,\tilde{A}\tilde{b})} \|\tilde{b} - \tilde{A} v\|^2 = \arg\min\limits_{\mathcal{K}_i(A,b)} \|y - X v\|^2, i=1,\dots,k.
\]
Thus it holds
\begin{equation*}
	\betaest_i = \arg\min\limits_{v \in \mathcal{K}_i(\tilde{A}^2,\tilde{A}\tilde{b})} \|\tilde{b} - \tilde{A} v\|^2,
\end{equation*}
Furthermore we have $\Sigma \beta(\eta_1) = \Sigma^{-1} P q$, i.e. the correct problem in the population is solved by $\beta(\eta_1)$ as well. Now we will restate the main result in \citet{Nemirovskii86} in our context:
\begin{theorem} Nemirovskii
	\label{th:Nemirovskii}
	\ \\
Assume that there are {$\tilde{\delta}=\tilde{\delta}(\nu,n)>0$, $\tilde{\epsilon}=\tilde{\epsilon}(\nu,n)>0$} such that for $\nu \in (0,1]$
	it holds that $
		\Prob\left(\|\Sigma - A^{1/2}\|_\calL \leq \tilde{\delta}\right) \geq 1- \nu/2, \;\;
		\Prob\left(\|\Sigma^{-1} Pq  - A^{-1/2}b\|  \leq \tilde{\epsilon}\right) \geq 1 - \nu/2$ and the conditions
	\begin{enumerate}
		\item 
		there is an $L=L(\nu,n)$ such that with probability at least $1-\nu/2$ it holds that
			$\max\left\{	\|A^{1/2}\|_\calL, \|\Sigma\|_\calL \right\} \leq L$,
		\item 
			there is a vector $u \in \R^k$ and constants $R, \mu>0$ such that
			$\beta(\eta_1) = \Sigma^{\mu} u$, $\|u\| \leq R$
	\end{enumerate}
	are satisfied. {If we stop according to the stopping rule $a^\ast$ as defined in (4) with $\tau \geq 1$ and $\zeta < \tau^{-1}$ then we have for any $\theta \in [0,1]$ with probability at least $1- \nu$}
	\[
		\left\|\Sigma^{\theta} \{\widehat{\beta}_{a^\ast} - \beta(\eta_1)\} \right\|^2
		\leq C^2(\mu,\tau, \zeta) R^{2(1-\theta)/(1+\mu)}\left(
			\tilde{\epsilon} + \tilde{\delta} R L^\mu
		\right)^{2(\theta+\mu)/(1+\mu)}.
	\]
\end{theorem}
{\it{Proof}}: Note first that on the set where $\|\Sigma - A^{1/2}\|_\calL \leq \tilde{\delta}$ holds with probability at least $1-\nu/2$ condition 1 also holds with $L=\|\Sigma\|_\calL + \tilde{\delta}$. Constrained on the set where all the conditions of the theorem hold with probability at least $1- \nu$ we consider Nemirovskii's $(\Sigma,A^{1/2},\Sigma^{-1} P q, A^{-1/2}b)$ problem with errors $\tilde{\delta}$ and $\tilde{\epsilon}$. Furthermore by assumption Nemirovskii's $(2\theta,R,L,1)$ conditions hold and thus the theorem follows by a simple application of the main theorem in \citet{Nemirovskii86}. \eop
\ \\
We will now apply Theorem \ref{th:Nemirovskii} to our problem. Due to the fact that $\eta_1>0$ it holds that $\Sigma^2$ is positive definite and thus invertible. We note that the spectral norm is dominated by the Frobenius norm. From Markov's inequality we get
\[
	\Prob\left(
		\|A - \Sigma^2 \| \geq \delta
	\right)  \leq \delta^{-2} \E\left(\left\|A -\Sigma^2\right\|^2\right).
\]
{Using Theorem 1, $\sum_{t=1}^n \|\tCov_i\|^4 \leq \|\tCov^2\|^2$ and
setting the right hand side to $\nu/2$ for $\nu \in (0,1]$ gives $\delta =  \nu^{-1/2} \|\tCov\|^{-2} \|\tCov^2\| C_\delta$}. In the same way $\epsilon =  \nu^{-1/2} \|\tCov\|^{-2}\|\tCov^2\| C_\epsilon$. Lemma \ref{lem:sqrRootConcentration} gives with probability at least $1-\nu/2$ the concentration results required by Theorem \ref{th:Nemirovskii} with
\begin{align*}
	\tilde{\delta} &=\nu^{-1/2} \frac{\|\tCov^2\|}{\|\tCov\|^2}C_\delta \{1 + o(1)\}\\
	\tilde{\epsilon} &=\left(\nu^{-1/2} \frac{\|\tCov^2\|}{\|\tCov\|^2}
		C_\epsilon + \nu^{-1}\frac{\|\tCov^2\|^2}{\|\tCov\|^4}C_\epsilon C_\delta\right)\{1+o(1)\}
\end{align*}
Conditions 1 and 2 of Theorem \ref{th:Nemirovskii} hold with a probability of at least $1-\nu/2$ by choosing $L = \tilde{\delta} + \|\Sigma\|_\calL$, $\mu=1$ and $R = \|\Sigma^{-3} P q\|$. Here we used that $\beta(\eta_1) = \Sigma^{-2} P q$. Thus the theorem yields for $\theta = 1$
\[
	\left\|\Sigma\{\beta(\eta_1) - \betaest_{a^\ast}\}\right\| \leq C(1,\tau,\zeta) \left(\tilde{\epsilon} + \tilde{\delta} R L \right).
\]
Denote $c(\tau,\zeta) = C(1,\tau,\zeta)\{1+o(1)\}$.
Finally we have $\|\Sigma^{-1}\|_\calL^{-1} \|v\| \leq \| \Sigma v\|$ for any $v \in \R^k$ and thus the theorem is proven with
\begin{align*}
 	c_1(\nu) &= \nu^{-1/2} c(\tau,\zeta) \|\Sigma^{-1}\|_\calL \left(
 		C_\epsilon + \|\Sigma\|_\calL\|\Sigma^{-3} P q\| C_\delta 
 	\right)\\
 	c_2(\nu) & = \nu^{-1}c(\tau,\zeta)\|\Sigma^{-1}\|_\calL \left(
 	C_\epsilon C_\delta + \|\Sigma^{-3} P q \| C_\delta^2
 	\right).
\end{align*} \eop

\subsection{Proof of Theorem 3}
The theorem is proved by contradiction.
Assume that $\widehat{\beta}_1 \longrightarrow
\beta_1$ in probability. Choosing $v \in \R^k$, $ v \neq 0$, orthogonal to $\beta_1$
implies that $v^\T \widehat{\beta}_1$ converges in probability to zero. Next we show that the second moment vanishes as well.

Let $M_d(z)=\max_{i \in \{1,\ldots,n\}^d}\E( \prod_{\nu = 1}^d z_{i_v}^2
	)$ for a random vector $z=(z_1,\ldots,z_n)^\T$ with existing mixed $(2d)$th moments. Using $(a+b)^2 \leq 2(a^2+b^2)$ for $a,b \in \R$ we obtain
\begin{align*}
	\E \left(v^\T b\right)^4 \leq& \frac{8^2 \|v\|^4}{\|\tCov\|^8} \E \left(
		\left\| 
			P N^\T \tCov^2 N q
		\right\|^4
		+
		\eta_1^4\left\|
			F^\T \tCov^2Nq
		\right\|^4
		+
		\eta_2^4\left\|
			PN^\T \tCov^2f
		\right\|^4
		+
		\eta_1^4\eta_2^4
		\left\|
			F^\T \tCov^2 f
		\right\|^4
	\right)\\
	\leq&
	8^2 \|v\|^4 \left\{M_4(N_1)  \|q\|^4 l^4 \|P\|^4 +
		M_2(N_1)M_2(F_1)\eta_1^4 \|q\|^4 l^2 k^2 \right.\\
		&\left.+ M_2(N_1)M_2(f_1)\eta_2^4 l^2 \|P\|^4 + M_2(F_1)M_2(f_1)\eta_1^4\eta_2^4 k^2  
	 \right\} < \infty, ~~  n \in \N.
\end{align*}
Thus, $\left(v^\T b\right)^2$ is uniformally integrable by the
theorem of de la Vall\'ee-Poussin and it follows that the directional
variance $\Var(v^\T b)$ has to vanish in the limit as well. 
Now, calculations similar to Theorem 1 yield
\begin{align*}
	\Var(v^\T b) =& \frac{\|\tCov^2\|^2}{\|\tCov\|^4} \left\{\eta_1^2 \|v\|^2\left(
		\|q\|^2 + \eta_2^2
	\right) + \|P^\T v\|^2 \left(
		\|q\|^2 + \eta_2^2
	\right) + (v^\T P q)^2\right\}\\
	\notag
	&+ \sum\limits_{t=1}^n\frac{\|\tCov_t\|^4}{\|\tCov\|^4} \sum\limits_{i=1}^l q_i^2\left(v^\T P_i\right)^2\{\E(N_{1,1}^4) - 3\},\;\;v \in \R^k.
\end{align*}

We assumed that $\|\tCov\|^{-2}\|\tCov^2\|$ does not converge to zero. It remains to check under which conditions $\Var(v^\T b)$ is larger than zero.
This will always be the case if $v \neq 0$ and $\eta_1 >
0$, $l=1$. For $\eta_1 = 0$ and $l>1$ a vector $v$ that lies in the range
of $P$ and is orthogonal to $\beta_1 \propto \spCC$ exists, thus
contradicting $\widehat{\beta}_1 \longrightarrow
\beta_1$ in probability.\eop

{\subsection{Proof of Theorem 4}
It is easy to verify that $\|\tCov\|^2 =
\mathrm{tr}(T^2) = n \gamma(0)$ and $
	\left\|\tCov^2\right\|^2 = n \gamma^2(0) + 2 \sum_{t=1}^{n-1} \gamma^2(t) (n-t) $.
If (\ref{eq:acfBound}) is fulfilled, then 
\begin{align*}
	n \gamma(0) \leq \left\|\tCov^2 \right\|^2 &\leq n \gamma^2(0)\left\{1 + 2 c^2\frac{1-\exp(-2\rho(n-1))}{\exp(2\rho)-1}\right\}
	\leq n \gamma^2(0)\left\{1 + \frac{2c}{\exp(2 \rho)-1}\right\}.
\end{align*}
It follows that $\|\tCov^2\| \sim n^{1/2}$.  \eop}

\subsection{Proof of Theorem 5}
Let $\gamma:\N \rightarrow \R$ be the autocovariance function of a stationary time series that has zero mean. For the autocovariance matrix $\tCov^2$ of the corresponding integrated process of order one we get
$
	\left[\tCov^2\right]_{t,s} = \sum_{i,j=1}^{t,s} \gamma(|i-j|),\;(t,s=1,\ldots,n).
$
Let $t \geq s$. By splitting the sum into parts with $i < j$ and $i >j$ we get
$
	\left[\tCov^2\right]_{t,s} = s \gamma(0) + \sum_{j=1}^s \sum_{i=1}^{t-j} \gamma(i) + \sum_{j=2}^s \sum_{i=1}^{j-1} \gamma(i).
$
Due to symmetry, $\left[\tCov^2\right]_{t,s} = \left[\tCov^2\right]_{s,t}$ for $s > t$.

First, consider the case that all $\gamma(j)$, $j>0$ are negative. Using (6) we obtain
\begin{align*}
	\notag
	\gamma(0) s \geq \left[\tCov^2\right]_{t,s} \geq \gamma(0)\left\{ s - c \sum\limits_{j=1}^s \sum\limits_{i=1}^{t-j} \exp(-\rho j) -
	c \sum\limits_{j=2}^s \sum\limits_{i=1}^{j-1} \exp(-\rho j)\right\}, ~ t\geq s.
\end{align*}
Evaluation of the geometric sums gives
\begin{align*}
	\notag
\left[\tCov^2\right]_{t,s} \geq \gamma(0)\left(s \left\{
		1 - \frac{2 c}{\exp(\rho) - 1}
	\right\} + c\frac{\exp(\rho)}{\left\{\exp(\rho) - 1 \right\}^2}
	\left\{
		 1 - \exp(-\rho s)
	\right\}\left[
		1 + \exp\{\rho (s- t)\}
	\right]\right).
\end{align*}
The second term on the right is always positive and the
positivity of the first term is ensured by the condition $\rho>\log(2c+1)$.
Hence, 
$
	\gamma(0)\left[1 - 2c\left\{\exp(\rho) -
    1\right\}^{-1}\right] s  \leq \left[\tCov^2\right]_{t,s} \leq \gamma(0)s, ~ s\geq 1.
$
If $\gamma(t)$, $t\geq 1$ is not purely negative, it can be bound by
\[
	\gamma(0)\left[1 - 2c\left\{\exp(\rho) -
    1\right\}^{-1}\right] s \leq 
    	\left[\tCov^2\right]_{t,s}
    \leq \gamma(0)
    \left[1 + 2c\left\{\exp(\rho) -
    1\right\}^{-1}\right] s.
\]
We write $\delta_1$ and $\delta_2$ for the constants in the lower and upper bound, respectively, so that $\delta_1 \min\{s,t\} \leq \left[\tCov^2\right]_{t,s} \leq \delta_2 \min\{s,t\}$ ($t,s=1,\ldots,n$).
This yields upper and lower bounds on the trace of $\tCov^2$ and shows that $\|\tCov\|^2 \sim n^2$.
Additionally, 
\begin{align*}
	\left[\tCov^4\right]_{t,t} &= \sum\limits_{l=1}^n \left[\tCov^2\right]_{t,l}\left[\tCov^2\right]_{l,t} =
	\sum\limits_{l=1}^t \left[\tCov^2\right]_{t,l}^2 + \sum\limits_{l=t+1}^n \left[\tCov^2\right]_{l,t}^2 \leq \frac{\delta_2^2}{6} t \left(6 n t - 4 t^2+3t + 1\right)\\
	\left[\tCov^4\right]_{t,t} &\geq \frac{\delta_1^2}{6} t \left(6 n t - 4 t^2+3t + 1\right).
\end{align*}
This implies upper and lower bounds on the trace of $\tCov^4$ in the form $c\, n(n+1)(n^2+n+1)$ for $c \in \{\delta_1^2/6,\delta_2^2/6\}$ and thus $\|\tCov^2\| \sim n^2$.
\eop

\subsection{Proof of Theorem 6}

First consider $n^{-1} X^\T \widehat{\tCov}^{-2} y$. Define $X_u =(X_{u,1},\ldots,X_{u,n})^\T = N P^\T + \eta_1 F$ and $y_u = (y_{u,1},\ldots,y_{u,n})^\T =N q + \eta_2 f$ such that $X = \tCov X_u$ and $y = \tCov y_u$. By the triangle inequality
\[
	\left\|n^{-1} X^\T \widehat{\tCov}^{-2} y - P q \right\| \leq \left\| n^{-1} X^\T \tCov^{-2} y - P q \right\| + \left\| n^{-1} X^\T \left( \widehat{\tCov}^{-2} - \tCov^{-2} \right) y \right\|.
\]
The first term on the right hand side is convergent to zero in probability due to Theorem 1.
The second term can be bounded as
\begin{align*}
	n^{-2}\left\| X^\T \left(\widehat{\tCov}^{-2}- \tCov^{-2}\right) y \right\|^2  \leq 
	 \|\tCov \widehat{\tCov}^{-2} \tCov - I_n\|_\calL^2 \, n^{-1} \|X_u^\T\|_\calL^2\, n^{-1} \|y_u\|^2.
\end{align*}
Since both $X_{u,1},\ldots,X_{u,n}$ and $y_{u,1},\ldots,y_{u,n}$ are independent  and identically distributed, it follows that $n^{-1} \|y_u\|^2$ is a strongly consistent estimator for $\E(y_{u,1}^2)$, as well as that $n^{-1} \| X_u^\T\|_\calL^2$ is bounded from above by $n^{-1} \|X_u^\T\|^2$, which is a strongly consistent estimator of $\E(\|X_{u,1}\|^2)$. Convergence in probability of $\left\|\tCov \widehat{\tCov}^{-2} \tCov - I_n\right\|^2_\calL$ to zero implies the convergence of $b(\widehat{\tCov})$ to $P q$ in probability.
To obtain the convergence rate $\| n^{-1} X^\T \tCov^{-2} y - P q \| = O_p(r_n)$, use Theorem 1 and $\|\tCov \widehat{\tCov}^{-2} \tCov - I_n\|_\calL =  O_p(r_n)$.         
The convergence of $\|n^{-1}X^\T \widehat{\tCov}^{-2} X - \Sigma^2\|$ is proven in a similar way.

To show the consistency and the rate of the corrected partial least squares estimator, we follow the same lines as in the proof of Theorem 2. First, $\delta =  r_n c_A(\nu)$ and $\epsilon = r_n  c_b(\nu)$ for $\nu \in (0,1]$ with constants $c_A(\nu)$, $c_b(\nu)$ are taken, such that
\begin{align*}
	\Prob\{\|A(\widehat{\tCov})^{1/2} - \Sigma\|_\calL &\leq r_n c_A(\nu)\} \geq 1 - \nu/2,\\ 
	\Prob\{\|A(\widehat{\tCov})^{-1/2} b(\widehat{\tCov}) - \Sigma^{-1} P q \| &\leq r_n c_b(\nu)\} \geq 1 - \nu/2.
\end{align*}
Moreover, $L = \|\Sigma\|_\calL + \delta$ and $R=\|\Sigma^{-3} P q\|$, $\mu=1$, satisfies conditions 1 and 2 in Theorem \ref{th:Nemirovskii} with probability at least $1- \nu/2$. Thus, with probability at least $1 - \nu$ we get by setting $\theta= 1$
\[
	\left\|
		\widehat{\beta}_{a^\ast}(\widehat{\tCov}) - \beta(\eta_1)
	\right\| \leq  r_n \, C(1, \tau, \zeta)\{1+o(1)\} \|\Sigma^{-1}\|_\calL \left[c_b(\nu) + c_A(\nu) \|\Sigma^{-3} P q\| \left\{\|\Sigma\|_\calL + r_n c_A(\nu)\right\}\right],
\]
where the constants $\zeta, \tau$ are taken from the definition of $a^\ast$.
\eop\

\section{Addendum to section 5, Simulations}
\begin{figure}
\captionsetup{width=0.6\textwidth}
	\figurebox{25pc}{25pc}{}[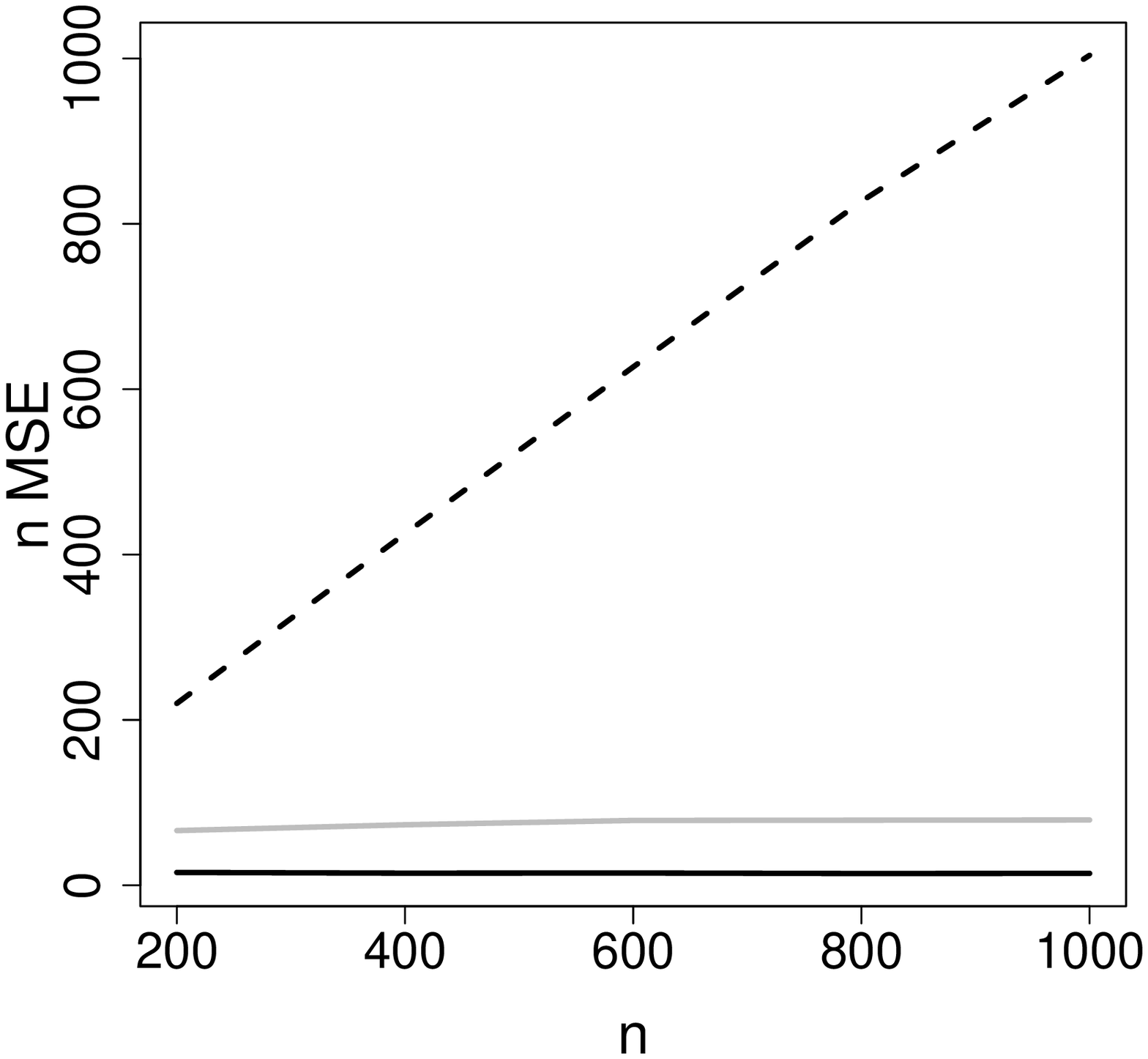]
   \caption[20pc]
   {Empirical mean squared eror of $\widehat{\beta}_1$ multiplied by $n$. The dependence structures are: autoregressive (grey), autoregressive integrated moving average (black, dashed) and corrected partial least squares on integrated data (black, solid). \label{fig:MSE}   
   }
\end{figure}

Figure \ref{fig:MSE} shows the differences in empirical mean squared error of $\hat\beta_1$ for various dependence structures considered in Section 5
in the setting with $l=i=1$. We calculated 
$$
n\mbox{MSE}(\hat\beta_1)=n\, {500}^{-1} \sum_{i=1}^{500}(\hat\beta_{1,i}-\beta_1)^2,
$$
where $\hat\beta_{1,i}$ denotes a partial least squares estimator in the $i$th Monte Carlo simulation based on $n$ observations.
If an
autoregressive dependence is present in the data and is ignored in the \PLS{}
algorithm, $n\mbox{MSE}(\hat\beta_1)$ is proportional to a
constant, which is larger than in the corrected \PLS{} case. Ignoring the integrated dependence in the data leads to
$n\mbox{MSE}(\hat\beta_1)$ growing linearly with $n$, which confirms our theoretical
findings in Section 3.   

\bibliography{bibliography}
\biblist